\documentclass[11pt]{article}
\usepackage{amssymb,amsfonts,amsmath}
\usepackage{epsfig}
\usepackage{graphicx}
\usepackage{color}
\usepackage{anysize}
\usepackage{hyperref}
\usepackage[french, english]{babel}

\marginsize{2.3cm}{2.3cm}{1cm}{2cm}

\begin{document}
\baselineskip=14pt

\makeatletter\@addtoreset{equation}{section}
\makeatother\def\theequation{\thesection.\arabic{equation}}

\newtheorem{defin}{Definition}[section]
\newtheorem{Prop}{Proposition}
\newtheorem{teo}{Theorem}[section]
\newtheorem{ml}{Main Lemma}
\newtheorem{con}{Conjecture}
\newtheorem{cond}{Condition}
\newtheorem{conj}{Conjecture}
\newtheorem{prop}[teo]{Proposition}
\newtheorem{lem}{Lemma}[section]
\newtheorem{rmk}[teo]{Remark}
\newtheorem{cor}{Corollary}[section]

\newcommand{\be}{\begin{equation}}
\newcommand{\ee}{\end{equation}}
\newcommand{\ben}{\begin{eqnarray}}
\newcommand{\benn}{\begin{eqnarray*}}
\newcommand{\een}{\end{eqnarray}}
\newcommand{\eenn}{\end{eqnarray*}}
\newcommand{\bp}{\begin{prop}}
\newcommand{\ep}{\end{prop}}
\newcommand{\bt}{\begin{teo}}
\newcommand{\et}{\end{teo}}
\newcommand{\bcor}{\begin{cor}}
\newcommand{\ecor}{\end{cor}}
\newcommand{\bcon}{\begin{con}}
\newcommand{\econ}{\end{con}}
\newcommand{\bcond}{\begin{cond}}
\newcommand{\econd}{\end{cond}}
\newcommand{\br}{\begin{rmk}}
\newcommand{\er}{\end{rmk}}
\newcommand{\bl}{\begin{lem}}
\newcommand{\el}{\end{lem}}
\newcommand{\bit}{\begin{itemize}}
\newcommand{\eit}{\end{itemize}}
\newcommand{\bd}{\begin{defin}}
\newcommand{\ed}{\end{defin}}
\newcommand{\bpr}{\begin{proof}}
\newcommand{\epr}{\end{proof}}

\newenvironment{proof}{\noindent {\em Proof}.\,\,}{\hspace*{\fill}$\halmos$\medskip}

\newcommand{\halmos}{\rule{1ex}{1.4ex}}
\def \qed {{\hspace*{\fill}$\halmos$\medskip}}

\newcommand{\fr}{\frac}
\newcommand{\Z}{{\mathbb Z}}
\newcommand{\R}{{\mathbb R}}
\newcommand{\E}{{\mathbb E}}
\newcommand{\C}{{\mathbb C}}
\renewcommand{\P}{{\mathbb P}}
\newcommand{\N}{{\mathbb N}}
\newcommand{\var}{{\mathbb V}}
\renewcommand{\S}{{\cal S}}
\newcommand{\T}{{\cal T}}
\newcommand{\W}{{\cal W}}
\newcommand{\X}{{\cal X}}
\newcommand{\Y}{{\cal Y}}
\newcommand{\h}{{\cal H}}
\newcommand{\f}{{\cal F}}

\renewcommand{\a}{\alpha}
\renewcommand{\b}{\beta}
\newcommand{\g}{\gamma}
\newcommand{\G}{\Gamma}
\renewcommand{\L}{\Lambda}
\renewcommand{\l}{\lambda}
\renewcommand{\d}{\delta}
\newcommand{\D}{\Delta}
\newcommand{\e}{\epsilon}
\newcommand{\s}{\sigma}
\newcommand{\B}{{\cal B}}
\renewcommand{\o}{\omega}
\newcommand{\eps}{\epsilon}

\newcommand{\nn}{\nonumber}
\renewcommand{\=}{&=&}
\renewcommand{\>}{&>&}
\newcommand{\<}{&<&}
\renewcommand{\le}{\leq}
\newcommand{\+}{&+&}

\newcommand{\pa}{\partial}
\newcommand{\ffrac}[2]{{\textstyle\frac{{#1}}{{#2}}}}
\newcommand{\dif}[1]{\ffrac{\partial}{\partial{#1}}}
\newcommand{\diff}[1]{\ffrac{\partial^2}{{\partial{#1}}^2}}
\newcommand{\difif}[2]{\ffrac{\partial^2}{\partial{#1}\partial{#2}}}

\title{Annealed vs Quenched Critical Points for a Random Walk Pinning Model}
\author{Matthias Birkner$^{\,1}$, Rongfeng Sun$^{\,2}$}
\date{Apr 14, 2009}
\maketitle

\footnotetext[1]{Dept.\ Biologie II, Abteilung Evolutionsbiologie,
University of Munich (LMU), Grosshaderner Str.\ 2, 82152 Planegg-Martinsried, Germany.
birkner@biologie.uni-muenchen.de}

\footnotetext[2]{Dept.\ Math., National University of Singapore, 2 Science Drive 2, Singapore, 117543. matsr@nus.edu.sg}

\begin{abstract}
We study a random walk pinning model, where conditioned on a simple random walk $Y$  on $\Z^d$
acting as a random medium, the path measure of a second independent simple random walk $X$ up
to time $t$ is Gibbs transformed with Hamiltonian $-L_t(X,Y)$, where $L_t(X,Y)$ is the
collision local time between $X$ and $Y$ up to time $t$. This model arises naturally
in various contexts, including the study of the parabolic Anderson model with moving catalysts,
the parabolic Anderson model with Brownian noise, and the directed polymer model. It
falls in the same framework as the pinning and copolymer models,
and exhibits a localization-delocalization transition as the inverse temperature $\beta$ varies.
We show that in dimensions $d=1,2$, the annealed and quenched critical values of $\beta$ are both 0,
while in dimensions $d\geq 4$, the quenched critical value of $\beta$ is strictly
larger than the annealed critical value (which is positive). This implies the existence of certain intermediate
regimes for the parabolic Anderson model with Brownian noise and the directed polymer model.
For $d\geq 5$, the same result has recently been established by Birkner, Greven and den Hollander
\cite{BGdH08} via a quenched large deviation principle. Our proof is based on a fractional moment
method used recently by Derrida, Giacomin, Lacoin and Toninelli \cite{DGLT07} to establish the
non-coincidence of annealed and quenched critical points for the pinning model in the disorder-relevant
regime. The critical case $d=3$ remains open.

\bigskip

\selectlanguage{french}
\begin{center}
{\bf R\'esum\'e}
\end{center}

Nous consid\'erons le mod\`ele de marche al\'eatoire avec {\it
  pinning} suivant : \'etant donn\'e une marche al\'eatoire simple $Y$
sur $\Z^{d}$ qui sert d'environnement al\'eatoire, on se donne une
mesure de Gibbs sur les trajectoires d'une marche al\'eatoire $X$
jusqu'au temps $t$ de Hamiltonien $-L_t(X,Y)$ o\`u $L_t(X,Y)$ est le
temps local d'intersection entre $X$ et $Y$ jusqu'au temps $t.$ Ce
mod\`ele appara\^it naturellement dans des contextes vari\'es tels que
l'\'etude du mod\`ele parabolique d'Anderson avec catalyseurs
mouvants, l'\'etude du mod\`ele parabolique d'Anderson avec bruit
Brownien ainsi que dans le cadre de l'\'etude de polym\`eres
dirig\'es. Ce mod\`ele appartient \`a la m\^eme classe que les
mod\`eles de {\it pinning} et copolym\`eres et pr\'esente une
transition localisation / d\'elocalisation quand la temp\'erature
inverse $\beta$ varie. Nous montrons qu'en dimension $d=1,2$ les
valeurs critiques {\it annealed} et {\it quenched } de $\beta$ sont
toutes deux 0 mais que en dimension $d\ge 4$ la valeur critique {\it
  quenched} de $\beta$ est strictement sup\'erieure \`a la valeur {\it
  annealed} (qui est positive). Ceci entraine l'existence de certains
r\'egimes interm\'ediaires pour le mod\`ele parabolique de Anderson
avec bruit Brownien et pour les polym\`eres dirig\'es. Pour $d\ge 5$
des r\'esultats similaires ont \'et\'e r\'ecemment \'etablis par
Birkner, Greven et den Hollander \cite{BGdH08} {\it via} un principe
de grandes d\'eviations {\it quenched}.  Notre preuve se fonde sur la
m\'ethode des moments fractionnaires utilis\'ee r\'ecemment par
Derrida, Giacomin, Lacoin et Toninelli \cite{DGLT07} pour \'etablir la
non-co\"incidence des valeurs critiques {\it quenched} et {\it
  annealed} du mod\`ele de pinning dans le r\'egime li\'e au
d\'esordre. Le cas de la dimension critique $d=3$ reste ouvert.
\bigskip

\selectlanguage{english}


\medskip


\end{abstract}

\section{Introduction and main result}

\subsection{The model and main results}
We first define the continuous time version of the {\it random walk pinning model}, which more precisely,
could be called the {\it random walk pinned to random walk model}. Let $X$ and $Y$ be two
independent continuous time simple random walks on $\Z^d$ with jump rates 1 and $\rho\geq0$ respectively.
Let $\mu_t$ denote the law of $(X_s)_{0\leq s\leq t}$. For $\beta \in\R$, which plays the role
of the inverse temperature (if $\beta>0$), and for a fixed realization of $Y$ acting as a random medium, we define
a Gibbs transformation of the path measure $\mu_t$. Namely, we define a new path measure
$\mu^\beta_{t,Y}$ on $(X_s)_{0\leq s\leq t}$ which is absolutely continuous w.r.t.\ $\mu_t$ with
Radon-Nikodym derivative
\be\label{muty}
\frac{{\rm d}\mu^\beta_{t,Y}}{{\rm d}\mu_t}(X) = \frac{e^{\beta L_t(X,Y)}}{Z^\beta_{t,Y}},
\ee
where $L_t(X,Y) = \int_0^t 1_{\{X_s=Y_s\}}{\rm d}s$ is the collision local time between $X$ and $Y$ up to
time $t$, and
\be\label{Zbetaty}
Z^\beta_{t,Y} = \E^X_0\big[e^{\beta L_t(X,Y)}\big]
\ee
is the {\it quenched partition function} which makes $\mu^\beta_{t,Y}$ a probability measure, where
$\E^X_x[\cdot]$ denotes expectation w.r.t.\ $X$ starting from $x\in\Z^d$. The
{\it quenched free energy} of the model is defined by
\be
F(\beta, \rho) = \lim_{t\to\infty} \frac{1}{t} \log Z^\beta_{t,Y}.
\ee
We will show below that the limit exists and is non-random. As a disordered system, it is also natural to consider
the {\it annealed partition function} $\E^Y_0[Z^\beta_{t,Y}]$ and the {\it annealed free energy}
\be
F_{\rm ann}(\beta, \rho) = \lim_{t\to\infty} \frac{1}{t} \log \E^Y_0[Z^\beta_{t,Y}].
\ee
Note that $\E^Y_0[Z^\beta_{t,Y}]=\E^{X-Y}_0[e^{\beta L_t(X-Y,0)}]$ is also the partition function of
a homogeneous pinning model (see e.g.\ Giacomin \cite{G07}), namely a random walk pinning model where the
random walk $X-Y$ (with jump rate $1+\rho$) is pinned to the site $0$ instead of to a random trajectory.

To define the discrete time version of the random walk pinning model, let $X,Y$ be discrete time simple
random walks on $\Z^d$. The Gibbs transformed path measure $\hat \mu^\beta_{N,Y}$, $N\in\N$, can be defined
similarly as in (\ref{muty}), where we replace $L_t(X,Y)$ by $L_N(X,Y)=\sum_{i=1}^N 1_{\{X_i=Y_i\}}$. We
then define $\hat Z^\beta_{N,Y}$, $\hat F(\beta)$, $\hat \mu^\beta_{N,{\rm ann}}$, $\hat F_{\rm ann}(\beta)$
similarly for the discrete time model as for the continuous time model. Note that the free energies
$\hat F(\beta)$ and $\hat F_{\rm ann}(\beta)$ now only depend on $\beta$ since there are no more
jump rates to adjust. To keep things simple, we focus only on $X$ and $Y$ being simple random walks in this
paper. However, we expect much of the same results to hold and the proofs to be adaptable for general random
walks, and we will comment on possible adaptations when appropriate.

Our first result is the existence of the quenched free energies $F(\beta, \rho)$ and $\hat F(\beta)$.
Existence of the annealed free energies $F_{\rm ann}(\beta,\rho)$ and $F_{\rm ann}(\beta)$ is well known
(see e.g.\ Chapter 2 in \cite{G07}). Before stating the result, we first introduce a two-parameter family of
constrained partition functions for the random walk pinning model, where apart from a shift in time for
the disorder $Y$, the random walk $X$ is subject to the constraint $X_t=Y_t$ in (\ref{muty}). In continuous
time setting, for $0<s<t<\infty$, define
\be\label{zpinned}
Z^{\beta, {\rm pin}}_{[s,t],Y} =
\E^X_{Y_s}\left[\exp\left\{\beta\int_0^{t-s} 1_{\{X_u=Y_{s+u}\}}du\right\} 1_{\{X_{t-s}=Y_t\}}\right].
\ee
For $0\leq m<n<\infty$ with $m,n\in\N_0$,
we define $\hat Z^{\beta, {\rm pin}}_{[m,n],Y}$ analogously for the discrete time model. For simplicity,
we will denote $Z^{\beta, {\rm pin}}_{[0,t],Y}$ by $Z^{\beta, {\rm pin}}_{t,Y}$, and
$\hat Z^{\beta, {\rm pin}}_{[0,N],Y}$ by $\hat Z^{\beta, {\rm pin}}_{N,Y}$.

\bt{\bf [Existence of quenched free energy]}\label{T:qexp}\\
For any $\beta\in\R$ and $\rho\geq 0$, there exists a non-random constant $F(\beta, \rho)$ such that
\be\label{lim1}
F(\beta, \rho) = \lim_{t\to\infty} \frac{1}{t} \log Z^\beta_{t,Y}
= \lim_{t\to\infty} \frac{1}{t}\log Z^{\beta, {\rm pin}}_{t,Y},
\ee
where the convergence are a.s.\ and in $L^1$ w.r.t.\ $Y$. Furthermore, we have the representation
\be\label{lim2}
F(\beta,\rho) = \sup_{t>0}\frac{1}{t} \, \E^Y_0\left[\log Z^{\beta, {\rm pin}}_{t,Y}\right].
\ee
Analogous statements hold for the discrete time model.
\et
\bcor{\bf [Existence of critical points]}\label{C:free}\\
There exist $0\leq \beta_{\rm c}^{\rm ann}\leq \beta_{\rm c}<\infty$ depending on $\rho\geq 0$ such that:
$F_{\rm ann}(\beta,\rho)=0$ if
$\beta<\beta^{\rm ann}_{\rm c}$ and $F_{\rm ann}(\beta, \rho)>0$ if $\beta>\beta^{\rm ann}_{\rm c}$;
$F(\beta,\rho)=0$ if $\beta <\beta_{\rm c}$ and $F(\beta,\rho)>0$ if $\beta>\beta_{\rm c}$. Analogous
statements hold for the discrete time model with annealed and quenched critical points
$\hat \beta^{\rm ann}_{\rm c}$ and $\hat \beta_{\rm c}$ respectively.
\ecor
{\bf Remark.} See (\ref{abccts}) and (\ref{abcdis}) for the exact values of $\beta^{\rm ann}_{\rm c}$ and
$\hat\beta^{\rm ann}_{\rm c}$.
\medskip

\noindent
{\bf Remark.} As in the pinning model (see e.g.\ \cite{G07}), $\beta_{\rm c}$ marks the transition between
a localized and a delocalized phase: when $\beta<\beta_{\rm c}$ and $F(\beta, \rho)=0$,
$L_t(X,Y)$ is typically of order $o(t)$ w.r.t.\ $\mu^\beta_{t,Y}$ for $t$ large; when $\beta>\beta_c$
and $F(\beta, \rho)>0$, $L_t(X,Y)$ is typically of order $t$ w.r.t.\ $\mu^\beta_{t,Y}$ for $t$ large.
Similarly, $\beta^{\rm ann}_{\rm c}$ marks the transition between the localized and delocalized phase
for the annealed homogeneous pinning model.

One question of fundamental interest in the study of disordered systems is to determine when is the
disorder strong enough to shift the critical point of the model, i.e., when is
$\beta^{\rm ann}_{\rm c} < \beta_{\rm c}$? For the pinning model, this question has recently been
essentially fully resolved independently by Derrida, Giacomin, Lacoin and Toninelli \cite{DGLT07},
and Alexander and Zygouras \cite{AZ08}. For the random walk pinning model, our main result is the
following.
\bt{\bf [Annealed vs quenched critical points]}\label{T:cpt}\\
In dimensions $d=1$ and $2$, we have $\beta^{\rm ann}_{\rm c}=\beta_{\rm c}=\hat \beta^{\rm ann}_{\rm c}=\hat\beta_{\rm c}=
0$. In dimensions $d\geq 4$, we have $0<\beta^{\rm ann}_{\rm c} <\beta_{\rm c}$ for each $\rho>0$ and
$0<\hat \beta^{\rm ann}_{\rm c} <\hat \beta_{\rm c}$. For $d\geq 5$, there exists $a>0$ s.t.\
$\beta_{\rm c}-\beta^{\rm ann}_{\rm c}\geq a\rho$ for all $\rho\in [0,1]$. For $d=4$ and for each $\delta>0$,
there exists $a_\delta>0$ s.t.\ $\beta_{\rm c}-\beta^{\rm ann}_{\rm c} \geq a_{\delta}\rho^{1+\delta}$ for all
$\rho\in [0,1]$.
\et

For purposes relevant to applications for the parabolic Anderson model with Brownian noise and the directed
polymer model, in $d\geq 4$, we prove instead a stronger version of Theorem \ref{T:cpt}. Define
\be\label{beta*}
\beta^*_{\rm c} = \sup\Big\{\beta \in\R : \sup_{t>0} Z^{\beta}_{t,Y} <\infty\ \mbox{a.s.  w.r.t. } Y\Big\}.
\ee
Define $\hat \beta^*_{\rm c}$ for the discrete time model analogously. Clearly $\beta^*_{\rm c}\leq \beta_{\rm c}$
and $\hat \beta^*_{\rm c} \leq \hat \beta_{\rm c}$. We have
\bt{\bf [Non-coincidence of critical points strengthened]}\label{T:cptst}\\
For $d\geq 4$,  we have $\beta^{\rm ann}_{\rm c}< \beta^*_{\rm c}$ for each $\rho>0$ and
$\hat \beta^{\rm ann}_{\rm c}<\hat \beta^*_{\rm c}$. For $d\geq 5$, there exists $a>0$ s.t.\
$\beta^*_{\rm c} - \beta^{\rm ann}_{\rm c} \geq a\rho$ for all $\rho\in [0,1]$. For $d=4$ and for
each $\delta>0$, there exists $a_\delta>0$ s.t.\ $\beta^*_{\rm c}-\beta^{\rm ann}_{\rm c}\geq a_\delta\rho^{1+\delta}$
for all $\rho\in [0,1]$.
\et
{\bf Remark.} Theorem \ref{T:cptst} for $d\geq 5$ (without bounds on the gap) has recently been established by
Birkner, Greven, and den Hollander \cite{BGdH08} as an application of a quenched large deviation principle for
renewal processes in random scenery. Our aim here is to give an alternative proof based on adaptations of the
fractional moment method used recently by Derrida et al \cite{DGLT07} in the pinning model context, and to extend
to the $d=4$ case. Loosely speaking, because $\P(X_n=Y_n) \sim Cn^{-d/2} = Cn^{-1-\alpha}$ by the local central
limit theorem, $d\geq 5$ corresponds to the case $\alpha>1$ in \cite{DGLT07}; $d=4$ corresponds to the case
$\alpha=1$, which was not covered in \cite{DGLT07}, but included in \cite{AZ08}; while $d=3$ corresponds to
the marginal case $\alpha=1/2$, which for the pinning model with Gaussian disorder was recently shown by Giacomin et
al~\cite{GLT08} to be disorder relevant. For the random walk pinning model, $d=3$ remains open.
\medskip

\noindent
{\bf Remark.} It is an interesting open question whether $\beta^*_{\rm c}=\beta_{\rm c}$, i.e., whether
the quenched partition function $Z^{\beta}_{t,Y}$ is uniformly bounded in $t$ a.s.\ w.r.t.\  $Y$ in the
entire delocalized phase. As communicated to us by F.L.Toninelli, this question also remains open for the
pinning and the copolymer models.
\medskip

Theorem \ref{T:cptst} for the continuous time model confirms Conjecture~1.8 of Greven and den Hollander \cite{GdH07}
(for $d\geq 4$)
that the parabolic Anderson model with Brownian noise could admit an equilibrium measure with an infinite
second moment. Theorem \ref{T:cptst} for the discrete time model can be used to disprove a conjecture of Garel and
Monthus \cite{GM06} that for the directed polymer model in random environment, the transition from weak to strong disorder
occurs at $\beta^{\rm ann}_{\rm c}$. See Sec.\ 1.4 for more details. For some special environments in special
dimensions, this conjecture has already been disproved by Camanes and Carmona \cite{CC07}. In Section 1.4, we
will show that the results of Derrida et al \cite{DGLT07} on the pinning model can also be used to disprove the
Garel-Monthus conjecture in $d\geq 4$. The reader can also consult Section 1.5 of Birkner et al~\cite{BGdH08} for more detailed
expositions on the implication of Theorem \ref{T:cptst} for the various models mentioned above.

In the remainder of the introduction, we point out a connection between the random walk pinning model and
the parabolic Anderson model with a single moving catalyst, and how does the random walk pinning model fit
in the same framework as the pinning and copolymer models. Lastly, we will introduce an {\it inhomogeneous random
walk pinning model} which generalizes both the pinning and the random walk pinning model.

\subsection{Parabolic Anderson model with a single moving catalyst}
As for the continuous time random walk pinning model, let $Y$ be a continuous time simple random walk on $\Z^d$
with jump rate $\rho\geq 0$. The parabolic Anderson model with a single moving catalyst is the solution of the following
Cauchy problem for the heat equation in a time-dependent random potential
\be \label{PAMs}
\begin{aligned}
\frac{\partial}{\partial t} u(t, x) &= \Delta u(t, x) + \beta \delta_{Y_t}(x)\, u(t, x),  \\
u(0,x) &= 1,
\end{aligned}
\qquad \qquad x\in \Z^d,\ t\geq 0,
\ee
where $\beta\in\R$ and $\Delta f(x) = \frac{1}{2d} \sum_{\Vert y-x\Vert =1} (f(y)-f(x))$ is the discrete
Laplacian on $\Z^d$. Heuristically, the time-dependent potential $\beta \delta_{Y_t}(x)$ can be interpreted
as a single catalyst with strength $\beta$ moving as $Y$, $u(t,x)$ is then simply the expected number of
particles alive at position $x$ at time $t$ for a branching particle system, where initially one particle
starts from each site of $\Z^d$, and independently, each particle moves on $\Z^d$ as a simple random walk,
and whenever the particle is at the same location as the catalyst $Y$, it splits into two particles with
rate $\beta$ if $\beta>0$ and is killed with rate $-\beta$ if $\beta<0$. For further motivations and a
survey on the parabolic Anderson model, see e.g.\ G\"artner and K\"onig \cite{GK05}.

Quantities of special interest in the study of the parabolic Anderson model are the quenched and annealed
$p$-th moment Lyapunov exponents.
\be
\lambda_0 = \lim_{t\to\infty}\frac{1}{t}\log u(t,0), \qquad \lambda_p=\lim_{t\to\infty}\frac{1}{t}\log \E^Y_0[u(t,0)^p].
\ee
The annealed $p$-th moment Lyapunov exponents for $p\in\N$ have been studied by G\"artner and Heydenreich in
\cite{GH06}. Here we show that
\bt{\bf [Existence of quenched Lyapunov exponent]}\label{T:qlyaexp}\\
For any $\beta\in\R$ and $\rho\geq 0$, there exists a non-random constant $\lambda_0=\lambda_0(\beta, \rho)$
such that for all $x\in\Z^d$,
\be\label{limlambda}
\lambda_0 = \lim_{t\to\infty} \frac{1}{t} \log u(t,x) \qquad \mbox{a.s.\ and in $L^1$ w.r.t. $Y$}.
\ee
Furthermore, $\lambda_0(\beta,\rho)=F(\beta,\rho)$, where $F(\beta,\rho)$ is as in (\ref{lim1}).
\et

Indeed, the solution of (\ref{PAMs}) admits the Feynman-Kac representation
\be\label{ufeynman}
u(t, x) = \E^X_x\left[\exp \left\{\beta \int_0^t 1_{\{X_{t-s}=Y_s\}}ds \right\}\right],
\ee
where $X$ is a simple random walk on $\Z^d$ with jump rate $1$ and $X_0=x$. Except for the time reversal of
$X$ in (\ref{ufeynman}), $u(t,x)$ has the same representation as that for $Z^\beta_{t,Y}$. The same proof
as for Theorem \ref{T:qexp} then applies, which gives rise to the same representation  for $\lambda_0$ as
for $F(\beta,\rho)$ in (\ref{lim2}) due to the fact that the variational expression in (\ref{lim2}) is
invariant w.r.t.\ time reversal for $X$.

\subsection{Relation to pinning and copolymer models}
We now explain in what sense does the random walk pinning model belong to the same framework as the pinning
and the copolymer models. For simplicity, we will examine the discrete time random walk pinning model with
a path measure associated with the partition function $\hat Z^{\beta, {\rm pin}}_{[0,N],Y}$, c.f.\ (\ref{zpinned}).

The pinning and copolymer models are both Gibbs transformation of a renewal process. More precisely,
let $\sigma=(\sigma_0=0, \sigma_1, \sigma_2, \cdots)$ be a renewal process on $\N_0$, where the inter-arrival times
$(\sigma_{i}-\sigma_{i-1})_{i\in\N}$ are i.i.d.\ $\N\cup\{\infty\}$-valued random variables with distribution
$\P(\sigma_1=i)=K(i)$ for some probability kernel $K$ on $\N\cup\{\infty\}$. Let $(\omega_i)_{i\in\N}$ be i.i.d.\
real-valued random variables with $\E[\omega_1]=0$ and $\E[e^{\lambda\omega_1}]<\infty$ for all $\lambda\in\R$.
Let $h\in\R$ and $\beta\geq0$. Then for a fixed $N\in\N$, the finite volume Gibbs weight for a given realization of
the renewal sequence $\sigma$ for both models are of the form
\be\label{Wform}
W(\sigma) = \left\{
\begin{aligned}
\prod_{i=1}^m w\big(\beta, h, (\omega_j)_{\sigma_{i-1}< j\leq \sigma_i}\big) \quad & \mbox{if } N=\sigma_m \mbox{\ for
  some\ } m\geq 1, \\
0 \qquad \qquad \qquad & \mbox{otherwise},
\end{aligned}
\right.
\ee
where
\be\label{115}
w\big(\beta, h, (\omega_j)_{0<j\leq n}\big) = \left\{
\begin{aligned}
e^{\beta\omega_n+h} \qquad \qquad \qquad & \qquad \mbox{pinning model}, \\
\frac{e^{\beta\sum_{j=1}^n(\omega_j+h)} + e^{-\beta\sum_{j=1}^n(\omega_j+h)}}{2} \quad &\qquad \mbox{copolymer model}.
\end{aligned}
\right.
\ee
See \cite{G07} for more on the pinning and copolymer models. For the discrete time random walk pinning model,
we can write
\be\label{hatzsplit}
\hat Z^{\beta, {\rm pin}}_{N,Y} = \E^X_0\big[e^{\beta L_N(X,Y)}1_{\{X_N=Y_N\}}\big]
= \sum_{m=1}^N \sum_{\sigma_0=0<\sigma_1<\cdots<\sigma_m=N}\prod_{i=1}^m \Big(e^\beta \P^X_0(\tau_{\theta_{\sigma_{i-1}}Y}=\sigma_i-\sigma_{i-1})\Big),
\ee
where $\theta_n Y=(Y_{n+i}-Y_n)_{i\in\N_0}$ denotes a shift in $Y$, and $\tau_Y=\tau_Y(X)= \min\{i\geq 1: X_i=Y_i\}$.
Let us denote $K(i) = \E^Y_0\big[\P^X_0(\tau_Y=i)\big]= \P^{X-Y}_0(\tau_0=i)$, then $K$ with
$K(\infty)=\P^{X-Y}_0(\tau_0=\infty)$ is the return time distribution of a renewal process on $\N_0$. Let $\Delta_i=Y_i-Y_{i-1}$.
We can then rewrite (\ref{hatzsplit}) as
\be\label{hatzsplit2}
\hat Z^{\beta, {\rm pin}}_{N,Y} = \sum_{m=1}^N \sum_{\sigma_0=0<\sigma_1<\cdots<\sigma_m=N}\prod_{i=1}^m
\Big(K(\sigma_i-\sigma_{i-1})\ w\big(\beta, (\Delta_j)_{\sigma_{i-1}<j\leq \sigma_i}\big)\Big),
\ee
where
\be\label{rwf}
w\big(\beta, (\Delta_i)_{0<i\leq n}\big) =\frac{e^\beta \P^X_0(\tau_Y=n)}{K(n)}, \qquad Y_i=\sum_{j=1}^i\Delta_j.
\ee
In view of (\ref{hatzsplit2}) and (\ref{rwf}), we see that the random walk pinning model associated with
$\hat Z^{\beta, {\rm pin}}_{[0,N],Y}$ is also a Gibbs transformation of a renewal process with inter-arrival
law $K$, except that the disorder $(\Delta_i)_{i\in\N}$ take values in $\Z^d$  and the Gibbs weight factor $w(\cdot)$
for each renewal gap has a more complicated dependence on the disorder than for the pinning and copolymer models.
Nevertheless, this simple observation motivates us to try to adapt the fractional moment method from the
pinning model to our context. In the actual proof, we will use an alternative representation for
$\hat Z^{\beta, {\rm pin}}_{[0,N],Y}$, as well as for $Z^{\beta, {\rm pin}}_{[0,t],Y}$, which
admits a simpler form for the weight factor $w(\cdot)$ than (\ref{rwf}). See (\ref{disz}) and (\ref{rwpartrep}).
We will see later on that despite the entirely different nature of the disorder, the random walk pinning model turns
out to be a close analogue of the pinning model. Lastly we note that the fractional moment method has recently
been successfully applied also to the copolymer model, see Bodineau, Giacomin, Lacoin and Toninelli \cite{BGLT08}
and Toninelli \cite{T08}.

\subsection{An inhomogeneous random walk pinning model}
Another common feature between the pinning and the random walk pinning model is that, for both models, the
annealed partition function is that of a homogeneous pinning model. A further intriguing interplay between
the two models is that we can define an {\it inhomogeneous random walk pinning model}, from which both models can be
obtained by partial annealing. More precisely, let $X$ and $Y$ be discrete time simple random walks on $\Z^d$,
let $(\omega_i)_{i\in\N}$ be i.i.d.\ real-valued random variables with $\E[\omega_1]=0$, and
$M(\lambda)=\log\E[e^{\lambda\omega_1}]$ is well-defined for all $\lambda\geq 0$. Let $h\in\R$
and $\beta\geq 0$. Then the discrete time inhomogeneous random walk pinning model is the Gibbs transformation of the
path measure $\mu_N$ of $X$ up to time $N$ with Radon-Nikodym derivative
\be\label{mixedpin}
\frac{{\rm d}\mu^{\beta,h}_{N,Y,\omega}}{{\rm d}\mu_N}(X) = \frac{\exp\big\{\sum_{i=1}^N (\beta\omega_i+h)1_{\{X_i=Y_i\}}\big\}}{Z^{\beta,h}_{N,Y,\omega}},
\ee
where $Z^{\beta,h}_{N,Y,\omega} = \E^X_0[\exp\big\{\sum_{i=1}^N (\beta\omega_i+h)1_{\{X_i=Y_i\}}\big\}]$ is the
partition function, and we now have two sources of disorder: the location of pinning as given by $Y$, and the strength
of pinning as given by $\beta\omega_i+h$. Note that under annealing w.r.t.\ $Y$,
\be\label{pinningeq}
\E^Y_0[Z^{\beta,h}_{N,Y,\omega}]= \E^{X-Y}_0\left[\exp\Big\{\sum_{i=1}^N(\beta\omega_i+h)1_{\{(X-Y)_i=0\}}\Big\}\right]
\ee
is the partition function of a pinning model (without boundary constraint $(X-Y)_N=0$), where the underlying
renewal process is given by the return times of $X-Y$ to $0$. On the other hand, under annealing w.r.t.\ $\omega$,
$$
\E^\omega[Z^{\beta,h}_{N,Y,\omega}] = \E^X_0\big[e^{(M(\beta)+h)L_N(X,Y)}\big]
$$
is the partition function of a random walk pinning model with parameter $M(\beta)+h$.

The continuous time version of the inhomogeneous random walk pinning model can be defined similarly with partition
function
$$
Z^{\beta,h}_{t,Y,B}=\E^X_0\Big[\exp\Big\{\beta\int_0^t 1_{\{X_s=Y_s\}}({\rm d}B_s+h{\rm d}s)\Big\}\Big],
$$
where $B_s$ is a standard Brownian motion.

The discrete time inhomogeneous random walk pinning model first appeared implicitly in Birkner \cite{B04}
in the study of the directed polymer model (the continuous time analogue can be found in Greven and den Hollander
\cite{GdH07}). Given a simple random walk $X$ on $\Z^d$, $\lambda \geq 0$, i.i.d.\ real-valued random variables
$(\omega(n,x))_{n\in\N, x\in\Z^d}$  with $M(\lambda')=\log\E[e^{\lambda'\omega(1,1)}]$
well-defined for all $\lambda'\geq 0$, the (normalized) partition function of the directed polymer model is given by
$$
Z^\lambda_{N,\omega} = \E^X_0\big[e^{\sum_{i=1}^N\{\lambda\omega(i,X_i)-M(\lambda)\}}\big].
$$
Note that $(Z^\lambda_{N,\omega})_{N\in\N}$ is a positive martingale. The critical point of the model can be defined by
$$
\lambda_c=\sup\{\lambda\geq 0: (Z^\lambda_{N,\omega})_{N\in\N} \text{ is uniformly integrable}\}
=\sup\{\lambda\geq 0: \lim\nolimits_{N\to\infty} Z^\lambda_{N,\omega} > 0 \ \text{a.s.}\}.
$$
In the literature, $[0,\lambda_c)$ and $(\lambda_c,\infty)$ are called respectively the weak and strong disorder
regimes, characterized respectively by the uniform integrability (or the lack of u.i.) of $(Z^\lambda_{N,\omega})_{N\in\N}$.
See \cite{CSY04} for an overview of the directed polymer model, and see \cite[Theorem 1.1 and Prop.~3.1]{CY06}
for the existence of $\lambda_c$. The Garel-Monthus conjecture \cite{GM06} asserts that $\lambda_c=\lambda_2:=
\sup\{\lambda\geq 0: \sup_{N\in\N}\E[(Z^\lambda_{N,\omega})^2]<\infty\}$.
On the other hand, Birkner \cite[Lemma~1]{B04} showed that if $Y$ is
an independent copy of $X$, and $(\tilde \omega(n,x))_{n\in\N,x\in\Z^d}$ is an i.i.d.\ field with
a tilted law
$\P(\tilde \omega(n,x) \in d\zeta) = e^{\lambda\zeta-M(\lambda)} \P(\omega(n,x) \in d\zeta)$,
independent of $X$, $Y$ and $\omega$, then the size-biased law of $Z^\lambda_{N,\omega}$ is the same
as the law of
\be\label{sizebias}
\tilde Z^{\lambda}_{N,\omega, \tilde \omega, Y} =
\E^X_0\Big[\exp\Big\{\sum_{i=1}^N \big(1_{\{X_i\neq Y_i\}} \lambda
\omega(i,X_i) + 1_{\{X_i=Y_i\}}\lambda \tilde\omega(i,X_i) -M(\lambda)\big)\Big\}\Big].
\ee
Namely, $\E[f(\tilde Z^\lambda_{N,\omega,\tilde\omega,Y})]=\E[Z^\lambda_{N,\omega}f(Z^\lambda_{N,\omega})]$
for all bounded $f:\R_+\to\R$. The uniform integrability of $(Z^\lambda_{N,\omega})_{N\in\N}$ is then
equivalent to the uniform tightness of the laws of $(Z^\lambda_{N,\omega, \tilde\omega, Y})_{N\in\N}$.
If we integrate out the disorder $\omega$ in (\ref{sizebias}), then
\be\label{panneal1}
\E[\tilde Z^\lambda_{N,\omega,\tilde \omega,Y} |\tilde \omega,Y] =
\E^X_0\big[e^{\sum_{i=1}^N (\lambda\tilde\omega(i,X_i)-M(\lambda))1_{\{X_i=Y_i\}}}\big]
\ee
is precisely the partition function of the inhomogeneous random walk pinning model. Further
integrating out $\tilde\omega$ gives the partition function of a random walk pinning model with
parameter $\hat\beta(\lambda)=M(2\lambda)-2M(\lambda)$,
$$
\E[\tilde Z^\lambda_{N,\omega,\tilde \omega,Y} |Y] = \E^X_0\big[e^{\sum_{i=1}^N (M(2\lambda)-2M(\lambda))1_{\{X_i=Y_i\}}}\big].
$$
Since $\E[(Z^\lambda_{N,\omega})^2]=\E[\tilde Z^\lambda_{N,\omega,\tilde\omega,Y}]$, $\hat\beta(\lambda_2)=\hat\beta^{\rm ann}_c$
with $\hat \beta^{\rm ann}_c$ being the annealed critical point as in Theorem \ref{T:cptst}. Since for
non-degenerate $\omega$, $\hat\beta(\lambda)$ is strictly increasing in $\lambda$, Theorem \ref{T:cptst} implies that in
$d\geq 4$, there exists $\lambda'>\lambda_2$ such that $\E[\tilde Z^{\lambda'}_{N,\omega,\tilde \omega,Y} |Y]$ is uniformly bounded
in $N$ a.s.\ w.r.t.\ $Y$. Therefore the law of $(\tilde Z^{\lambda'}_{N,\omega,\tilde\omega,Y})_{N\in\N}$ is uniformly
tight, and hence $\lambda_c\geq \lambda'>\lambda_2$, which disproves the conjecture of Garel and Monthus~\cite{GM06}.
Since our proof is based on bounding fractional moments, we will in fact exhibit a $\lambda'>\lambda_2$ such that
$$
\sup_{N\in\N} \E\Big[\E[\tilde Z^{\lambda'}_{N,\omega,\tilde \omega,Y} |Y]^\gamma\Big] <\infty \qquad
\mbox{for some } \gamma \in (0,1).
$$
See (\ref{155}). Since $\tilde Z^\lambda_{N,\omega,\tilde \omega,Y}$ is the size-biased version of the partition function
$Z^{\lambda'}_{N,\omega}$ of the directed polymer model,
$$
\E\big[(Z^{\lambda'}_{N,\omega})^{1+\gamma}\big] = \E\big[\big(\tilde Z^{\lambda'}_{N,\omega,\tilde \omega,Y}\big)^\gamma\big]
\leq \E\big[\E[\tilde Z^{\lambda'}_{N,\omega,\tilde \omega,Y} |Y]^\gamma\big].
$$
Therefore, beyond the regime of $\lambda$ where $Z^\lambda_{N,\omega}$ is a $L_2$ bounded martingale, there is a regime where
$Z^\lambda_{N,\omega}$ has uniformly bounded $(1+\gamma)$-th moment for some $\gamma\in (0,1)$.

Finally, we point out that based on (\ref{sizebias}), the results of Derrida et al~\cite{DGLT07} for the pinning model
can also be used to disprove the Garel-Monthus conjecture in $d\geq 4$:
In (\ref{panneal1}), conditioned on $Y$, $(\tilde\omega(i,Y_i))_{1\leq i\leq N}$ are i.i.d.
Therefore if we fix an i.i.d.\ sequence $(\bar\omega_i)_{i\in\N}$ equally distributed with $\tilde\omega(1,1)$,
then $\E[\tilde Z^\lambda_{N,\omega,\tilde\omega,Y}|\tilde\omega,Y]$ is equally distributed with
$$
\E^X_0[e^{\sum_{i=1}^N (\lambda\bar\omega_i-M(\lambda))1_{\{X_i=Y_i\}}}].
$$
Integrating out $Y$ then gives the partition of a pinning model,
\be\label{pinpartition}
Z^{\beta,h}_{N,\bar\omega}=\E^{X-Y}_0[e^{\sum_{i=1}^N (\lambda\bar\omega_i-M(\lambda))1_{\{(X-Y)_i=0\}}}]
\ee
with parameters $\beta(\lambda)=\lambda$, $h(\lambda)=-M(\lambda)$ (c.f.\ (\ref{pinningeq})), and underlying renewal process
$K(n)=\P^{X-Y}_0(\tau_0=n)$ where $\tau_0$ is the first return time of $X-Y$ to $0$. It is easy to check
that the critical curve for the annealed pinning model is given by
$h^{\rm ann}_c(\beta)=M(\lambda)-M(\lambda+\beta)-\log\P^{X-Y}_0(\tau_0<\infty)$. By the definition of $\lambda_2$, $(\beta(\lambda_2),h(\lambda_2))$
lies on this annealed critical curve. Since in $d\geq 4$, $K(n)\sim c n^{-\frac{d}{2}}$ has tail exponent
$\alpha =\frac{d}{2}-1 \geq 1$, it follows from Derrida et al \cite{DGLT07} that there exists a continuous
curve $h^*(\beta)$ strictly above $h^{\rm ann}_c(\beta)$, such that for all $h\leq h^*(\beta)$,
$\E[(Z^{\beta,h}_{N,\bar\omega})^\gamma]$ is uniformly bounded in $N$ for some $\gamma \in (0,1)$. Therefore we can choose $\lambda'>\lambda_2$ such that
$-M(\lambda')\leq h^*(\lambda')$, and hence $\E[(Z^{\lambda', -M(\lambda')}_{N,\bar\omega})^\gamma]$ is uniformly bounded in $N$ for some
$\gamma \in (0,1)$. By the same reasoning as before, this implies the uniform tightness of
$(\tilde Z^{\lambda'}_{N,\omega,\tilde\omega,Y})_{N\in\N}$, and hence
$\lambda_c\geq \lambda'>\lambda_2$. We remark that in \cite{DGLT07}, only the constrained version of the partition function
$Z^{\beta,h}_{N,\bar\omega}$ is considered, i.e., the constraint $1_{\{X_N=Y_N\}}$ is inserted in
(\ref{pinpartition}). However, the proof there can be easily adapted to the non-constrained version, as can be seen later in
our analysis of the random walk pinning model.
Most recently, Giacomin, Lacoin and Toninelli~\cite{GLT08} extended their fractional moment technique to the pinning model
with Gaussian disorder in the critical dimension, i.e., $K(n)\sim c n^{-\frac{3}{2}}$, which corresponds to $d=3$ for the random walk pinning
model considered here. Except for the technical point that \cite{GLT08} only considered the constrained pinning model, their result would
imply $\lambda_c>\lambda_2$ for the directed polymer model in Gaussian environment in $d=3$, since in (\ref{pinpartition}), the
exponentially tilted law of a Gaussian is a shifted Gaussian.

\subsection{Outline}
The rest of the paper is organized as follows. In Section 2, we prove Theorem \ref{T:qexp}, Corollary \ref{C:free},
and Theorem \ref{T:qlyaexp}. In Section 3, we prove Theorem \ref{T:cpt} for $d=1,2$. In Section 4, we prove Theorem
\ref{T:cptst} in the discrete time case. Lastly in Section 5, we prove Theorem \ref{T:cptst} in the continuous time
case. The proof of Theorem \ref{T:cptst} does not rely on the existence of the quenched free energies. Readers
interested in how the fractional moment method is applied in this context can go directly to Sections 4 and 5.

\section{Existence of the quenched free energy}
In this section, we prove Theorems \ref{T:qexp}, \ref{T:qlyaexp} and Corollary \ref{C:free}.
\medskip

\noindent
{\bf Proof of Theorem \ref{T:qexp}.}
We consider first the constrained partition functions $\hat Z^{\beta, \rm pin}_{N,Y}$ and $Z^{\beta,\rm pin}_{t,Y}$.
For the discrete time model, by the super-additive ergodic theorem (see e.g.\ Sec.\ 6.6 in Durrett \cite{D96}) applied to
$(\log \hat Z^{\beta, \rm pin}_{[m,n],Y})_{0\leq m< n}$, we have
$$
\hat F(\beta) = \lim_{n\to\infty} \frac{1}{n} \log \hat Z^{\beta,\rm pin}_{n,Y} = \sup_{n\in\N} \E^Y_0\big[\log \hat Z^{\beta, \rm pin}_{n,Y}\big],
$$
where the convergence is a.s.\ and in $L^1$. For the continuous time model, we have to apply the super-additive ergodic theorem first
along the integer times, and then extend the convergence along all real times. Clearly $(\log Z^{\beta, {\rm pin}}_{[m,n],Y})_{0\leq m<n}$
satisfies all the conditions of the super-additive ergodic theorem, therefore
\be\label{lamconv}
F^{\rm pin}(\beta) = \lim_{n\to\infty} \frac{1}{n} \log Z^{\beta, {\rm pin}}_{n,Y} =
\sup_{n\in\N} \frac{1}{n} \E^Y_0\big[\log Z^{\beta,{\rm pin}}_{n,Y}\big] \qquad a.s.\ \mbox{ and in } L^1.
\ee
To extend the a.s.\ convergence to real $t\to\infty$, we need the following crude estimates.

\bp\label{P:LIL}
Let $(X_t)_{t\geq 0}$ be a continuous time random walk on $\Z^d$ with jump rate $1$. Let $\Vert\cdot\Vert_1$
denote $L^1$ norm in $\Z^d$. Then
\begin{itemize}
\item[$(i)$] There exists $C>0$ such that a.s.\ $\Vert X_t\Vert_1 < C\sqrt{t\log\log t}$ for all $t$ sufficiently large.
\item[$(ii)$] $\P^X_0(X_t =x) \geq C(1+t)^{-\frac{d}{2}}(2d)^{-\Vert x\Vert_1}$ uniformly for all $t>0$
and $x\in\Z^d$ with $\Vert x\Vert_1 \leq t/2$.
\end{itemize}
\ep
{\bf Proof.} Part (i) is a consequence of the law of the iterated logarithm. Part (ii) follows by forcing $X$
to visit $x$ after exactly $\Vert x\Vert_1$ number of jumps, and then return to $x$ at time $t$. The factor
$(1+t)^{-\frac{d}{2}}$ arises from the local central limit theorem.
\qed

Note that for $t\geq 1$, by super-additivity, we have
\be
\frac{1}{t}\left(\log Z^{\beta, {\rm pin}}_{\lfloor t-t^{2/3}\rfloor, Y}
+ \log Z^{\beta, {\rm pin}}_{[\lfloor t-t^{2/3}\rfloor, t], Y}\right)
\leq \frac{1}{t} \log Z^{\beta, {\rm pin}}_{t,Y}
\leq \frac{1}{t}\left(\log Z^{\beta, {\rm pin}}_{\lfloor t+t^{2/3}\rfloor, Y} -
\log Z^{\beta, {\rm pin}}_{[t, \lfloor t+t^{2/3}\rfloor], Y}\right).
\ee
By (\ref{lamconv}), a.s.\
$F^{\rm pin} = \displaystyle \lim_{t\to\infty} t^{-1} \log Z^{\beta, {\rm pin}}_{\lfloor t-t^{2/3}\rfloor, Y}
= \lim_{t\to\infty} t^{-1} \log Z^{\beta, {\rm pin}}_{\lfloor t+t^{2/3}\rfloor, Y}$. On the other hand,
\be
\begin{aligned}
Z^{\beta, {\rm pin}}_{[\lfloor t-t^{2/3}\rfloor, t],Y} &\ \leq\
e^{|\beta| (t- \lfloor t-t^{2/3}\rfloor)} \ \P^X_0\big(X_{t-\lfloor t-t^{2/3}\rfloor}=Y_t - Y_{\lfloor t-t^{2/3}\rfloor}\big),\\
Z^{\beta, {\rm pin}}_{[\lfloor t-t^{2/3}\rfloor, t],Y} &\ \geq\
e^{-|\beta| (t- \lfloor t-t^{2/3}\rfloor)} \ \P^X_0\big(X_{t-\lfloor t-t^{2/3}\rfloor}=Y_t - Y_{\lfloor t-t^{2/3}\rfloor}\big).
\end{aligned}
\ee
By Proposition \ref{P:LIL}, for $t$ sufficiently large, $\Vert Y_t - Y_{\lfloor t-t^{2/3}\rfloor}\Vert_1 \leq 2C\sqrt{t
  \log\log t} < \frac{(t-\lfloor t-t^{2/3}\rfloor)}{2}$, and hence
\be\nonumber
\P^X_0\big(X_{t-\lfloor t-t^{2/3}\rfloor}=Y_t - Y_{\lfloor t-t^{2/3}\rfloor}\big)
\geq C \big(1+t-\lfloor t-t^{2/3}\rfloor\big)^{-d/2} (2d)^{-2C\sqrt{t\log\log t}},
\ee
from which we obtain
$\displaystyle \lim_{t\to\infty} t^{-1} |\log Z^{\beta, {\rm pin}}_{[\lfloor t-t^{2/3}\rfloor, t],Y}|=0$.
Similarly,
$\displaystyle \lim_{t\to\infty} t^{-1} |\log Z^{\beta, {\rm pin}}_{[t, \lfloor t+t^{2/3}\rfloor],Y}|=0$.
This establishes the a.s.\ convergence in (\ref{lamconv}) for $t\to\infty$ in place of $n\to\infty$ for $n\in\N$.
To obtain $L^1$ convergence, it remains to verify the uniform integrability of $(t^{-1}\log Z^{\beta, {\rm pin}}_{t,Y})_{t\geq 1}$.
Note that
$$
t^{-1}\log Z^{\beta, {\rm pin}}_{t,Y} \leq t^{-1} \log Z^\beta_{t,Y} \leq \beta \qquad a.s.\ w.r.t.\ Y,
$$
while
$$
t^{-1} \log Z^{\beta, \rm pin}_{t,Y} \geq t^{-1} \log p_t(Y_t),
$$
where $p_t$ denotes the transition kernel of $X$. Using estimates (\ref{pest2})--(\ref{pest4}) below, it is easy
to see that $|t^{-1} \log p_t(Y_t)|_{t\geq 1}$ is uniformly integrable, hence $(t^{-1}\log Z^{\beta, {\rm pin}}_{t,Y})_{t\geq 1}$
is also uniformly integrable. Note that because $\log Z^\beta_{t,Y}\geq 0$, the unconstrained partition function
$(t^{-1}\log Z^\beta_{t,Y})_{t\geq 0}$ is also uniformly integrable.

We now consider the unconstrained partition functions $\hat Z^\beta_{N,Y}$ and $Z^\beta_{t,Y}$. The argument is the same
for discrete and continuous times, so we only consider the latter. Clearly
$Z^{\beta}_{t,Y} \! >\! Z^{\beta, {\rm pin}}_{t,Y}$. To upper bound $Z^{\beta}_{t,Y}$ in terms of $Z^{\beta, \rm pin}_{t,Y}$,
we can let $X$ run freely until time $t-t^{3/4}$ ($3/4$ is somewhat ad hoc), which gives a contribution of order
$Z^\beta_{t-t^{3/4},Y}$, and then force $X$ to go to $Y_t$ at time $t$. If $X_{t-t^{3/4}}$ is not too far from $Y_{t-t^{3/4}}$,
then we expect the cost of forcing $X_t=Y_t$ to be negligible, and if such $X$ gives the dominant contribution in $Z^\beta_{t-t^{3/4},Y}$,
then we are essentially done.

We now make the above heuristics precise. Note that
\be\label{112}
Z^{\beta}_{t,Y} \leq e^{|\beta| t^{3/4}} Z^\beta_{t-t^{3/4},Y}.
\ee
We claim that for $t$ sufficiently large,
\be\label{113}
\E^X_0\left[e^{\beta L_{t-t^{3/4}}(X,Y)} 1_{\{\Vert X_{t-t^{3/4}}\Vert_1\leq t^{2/3}\}}\right]
\geq \E^X_0\left[e^{\beta L_{t-t^{3/4}}(X,Y)} 1_{\{\Vert X_{t-t^{3/4}}\Vert_1 > t^{2/3}\}}\right].
\ee
By Proposition \ref{P:LIL}, for $t$ sufficiently large, we have
$\sup_{0\leq s\leq t} \Vert Y_s\Vert_1 \leq C\sqrt{t\log\log t}$.
Define recursively stopping times $\sigma_1=0$, and for $n\in\N$,
\be
\begin{aligned}
\tau_n &= \inf \{s \in (\sigma_n, t- t^{3/4}] : \Vert X_s\Vert_1 \geq t^{2/3}/2\},  \\
\sigma_{n+1} &= \inf \{s \in (\tau_n, t-t^{3/4}] : \Vert X_s \Vert_1 \leq C\sqrt{t\log\log t}\},
\end{aligned}
\ee
where we set $\sigma_n$, $\tau_n$ to $t-t^{3/4}$ if the infimum is taken over an empty set. Then
\begin{eqnarray}
&& \!\!\!\!\!\!\!\!\!\!\!\!\!\!\!\!\!\!
\E^X_0\left[e^{\beta L_{t-t^{3/4}}(X, Y)}1_{\{\Vert X_{t-t^{3/4}}\Vert_1 > t^{2/3}\}}\right]
=
\sum_{n=1}^\infty \E^X_0\left[e^{\beta L_{\tau_n}(X,Y)}1_{\{\tau_n
    <\sigma_{n+1}=t-t^{3/4}, \Vert X_{t-t^{3/4}}\Vert_1 > t^{2/3}\}}\right] \nonumber \\
&& \quad \quad =\quad
\sum_{n=1}^\infty \E^X_0\left[e^{\beta L_{\tau_n}(X,Y)} 1_{\{\tau_n <t-t^{3/4}\}}\,
\P^X_0\big(\sigma_{n+1} =t-t^{3/4}, \Vert X_{t-t^{3/4}}\Vert_1 >
  t^{2/3} \big| X_{\tau_n}\big)\right] \nonumber \\
&&\quad \quad\leq\quad
\sum_{n=1}^\infty \E^X_0\left[e^{\beta L_{\tau_n}(X, Y)}1_{\{\tau_n <t-t^{3/4}\}}\,
\P^X_0\big(\sigma_{n+1} =t-t^{3/4}, \Vert X_{t-t^{3/4}}\Vert_1 \leq
  t^{2/3}\big| X_{\tau_n}\big)\right] \nonumber\\
&&\quad \quad \leq\quad
\E^X_0\left[e^{\beta L_{t-t^{3/4}}(X, Y)}1_{\{\Vert X_{t-t^{3/4}}\Vert_1 \leq t^{2/3}\}}\right],
\end{eqnarray}
where in the first inequality we used the fact that $t^{2/3}/2>> \sqrt{t\log\log t} >> \sqrt{t}$ for $t$ large.
This proves the claim (\ref{113}). By Proposition \ref{P:LIL}, we have
$\P^X_0(X_t=Y_t | X_{t-t^{3/4}}=x)\geq C(1+t^{3/4})^{-d/2} (2d)^{-2t^{2/3}}$ uniformly for
$\Vert x\Vert_1 \leq t^{2/3}$. Hence
\be\nonumber
Z^{\beta, {\rm pin}}_{t,Y}\ \geq\
C (1+t^{3/4})^{-d/2}(2d)^{-2t^{2/3}} e^{-|\beta|t^{3/4}}
\E^X_0\left[e^{\beta L_{t-t^{3/4}}(X,Y)}1_{\{\Vert X_{t-t^{3/4}}\Vert_1 \leq t^{2/3}\}}\right].
\ee
Combined with (\ref{112}) and (\ref{113}), we find
\be\nonumber
Z^\beta_{t,Y} \leq 2 C^{-1} (1+t^{3/4})^{d/2}(2d)^{2t^{2/3}}e^{2|\beta| t^{3/4}} Z^{\beta, {\rm pin}}_{t,Y}.
\ee
Since $Z^{\beta}_{t,Y}>Z^{\beta, {\rm pin}}_{t,Y}$, (\ref{lim1}) follows with $F(\beta,\rho)=F^{\rm pin}$.

Lastly, (\ref{lim2}) holds because (\ref{lamconv}) is valid with $F^{\rm pin}=F(\beta,\rho)$ if we take the
limit in (\ref{lamconv}) along $nt$, $n\in\N$, for any fixed $t>0$.
\qed

\bigskip
\noindent
{\bf Proof of Corollary \ref{C:free}.} From the theory for homogeneous pinning models (see e.g.\ Chapter 2 of
\cite{G07}), it is known that $\beta^{\rm ann}_{\rm c}$ exists, and $\beta^{\rm ann}_{\rm c}=0$ if the renewal process
underlying the pinning model is recurrent (i.e., the random walk $X-Y$ is recurrent), and
$\beta^{\rm ann}_{\rm c}>0$ if the random walk $X-Y$ is transient. The statement
$\beta^{\rm ann}_{\rm c}\leq \beta_{\rm c}$ follows from
$$
F(\beta,\rho) = \lim_{t\to\infty} t^{-1} \E^Y_0\big[\log Z^\beta_{t,Y}\big]
\leq \lim_{t\to\infty} t^{-1} \log \E^Y_0[Z^\beta_{t,Y}] = F_{\rm ann}(\beta, \rho)
$$
by the $L^1$ convergence in Theorem \ref{T:qexp} and Jensen's inequality. The statement $\beta_c\geq 0$
follows from the fact that for $\beta<0$, $F(\beta, \rho)=0$. Indeed, for $\beta<0$, $Z^\beta_{t,Y} \leq 1$, while
$$
\log Z^\beta_{t,Y} =\log\E^X_0\big[e^{\beta L_t(X,Y)}\big] \geq \beta \E^X_0[L_t(X,Y)]\geq \beta \int_0^t \frac{C}{(1+s)^{d/2}}{\rm d}s =o(t),
$$
where we used the local central limit theorem that $\P^X_0(X_t=x) \leq C(1+t)^{-d/2}$ uniformly in $t>0$ and $x\in\Z^d$.
The existence and finiteness of $\beta_{\rm c}$ then follows from (\ref{lim2}) and the monotonicity of $F(\beta,\rho)$ in $\beta$.
The proof for the discrete time model is identical.
\qed

\bigskip
\noindent
{\bf Proof of Theorem \ref{T:qlyaexp}.} The difference between the Feynman-Kac representation of
$u(t,x)$ in (\ref{ufeynman}) and the representation for $Z^{\beta}_{t,Y}$ in (\ref{Zbetaty}) is: (1)
time-reversal for $X$; (2) in (\ref{ufeynman}), $X$ starts at $x$ instead of on $Y$. The same
proof as for Theorem \ref{T:qexp} shows that
$\displaystyle \lim_{t\to\infty} t^{-1} u(t, Y_t)=F(\beta,\rho)$ a.s.\ w.r.t.\ $Y$ where $F(\beta,\rho)$
is as in (\ref{lim2}). To compare $u(t,x)$ with $u(t,Y_t)$, note that
\be
u(t,x) \geq \P^X_0\big(X_{t^{2/3}}=Y_{t-t^{2/3}}-x\big) e^{-|\beta|t^{2/3}} u\big(t-t^{2/3}, Y_{t-t^{2/3}}\big),
\ee
which a.s.\ gives the correct lower bound on the exponential scale as $t\to\infty$. For the upper bound, note that
if $\beta\leq 0$, then $u(t,x)\leq 1$, which suffices by Corollary \ref{C:free}. If $\beta>0$, then for any $\eps>0$,
a.s.\ we can find $T_{\eps,Y}$ sufficiently large s.t.\ $t^{-1}\log u(t,Y_t) \leq F(\beta,\rho) +\eps$ for all
$t\geq T_{\eps,Y}$. In (\ref{ufeynman}), let $\tau=\inf\{s\in [0,t]: X_s=Y_{t-s}\}$ with $\tau=t$ if the set is
empty. Then for all $t>T_{\eps,Y}$ and $x\in\Z^d$, we have
\be
u(t,x) \leq \P^X_x(\tau \geq t-T_{\eps,Y}) e^{\beta T_{\eps,Y}}+ \E^X_x\big[u(t-\tau, Y_{t-\tau}) 1_{\{\tau<t-T_{\eps,Y}\}}\big]
\leq e^{\beta T_{\eps,Y}} + e^{(F(\beta,\rho)+\eps)t}.
\ee
Since $\eps>0$ can be arbitrarily small, a.s.\ this provides the correct upper bound for $u(t,x)$ on the exponential scale
as $t\to\infty$. The $L^1$ convergence in (\ref{limlambda}) follows from the uniform boundedness of $|\log u(t,x)|$ in
$t$, $x$ and $Y$.
\qed

\section{Coincidence of critical points in $d=1$ and $2$}
{\bf Proof of Theorem \ref{T:cpt} for $d=1$ and $2$.} The proof for the discrete and continuous time cases are essentially
the same, except that the estimates for the continuous time random walk transition kernel is slightly more involved. So we will only
consider the continuous time case. As pointed out in the proof of Corollary \ref{C:free}, because the random walk $X-Y$ is recurrent
in $d=1$ and $2$, $\beta^{\rm ann}_{\rm c}=0$. By (\ref{lim2}), to show $\beta_{\rm c}=0$, it suffices to show that
for any $\beta>0$, there exists $t>0$ such that
$\E^Y_0[\log Z^{\beta, {\rm pin}}_{t,Y}]>0$. We can write
\be \label{125}
\E^Y_0\big[\log Z^{\beta, {\rm pin}}_{t,Y}\big] =
\E^Y_0\big[\log \P^X_0(X_t=Y_t)\big] + \E^Y_0\left[\log \E^X_0\big[e^{\beta L_t(X,Y)} \big | X_t=Y_t\big]\right].
\ee
We first we estimate  $\E^Y_0\left[\log \P^X_0(X_t=Y_t)\right] =\sum_{x\in\Z^d} p_{\rho t}(x) \log p_t(x)$,
where $p_t(x)$ denotes the transition probability of a jump rate 1 continuous time simple random walk on $\Z^d$.
We then find lower bounds for the second term in (\ref{125}) for $d=1$ and $d=2$.

\bl\label{L:srwest}
For all $\rho\geq 0$, we have
\be\label{126}
\lim_{t\to\infty} \frac{\sum_{x\in\Z^d} p_{\rho t}(x) \log p_{t}(x)}{\log t} = -\frac{d}{2}.
\ee
\el
{\bf Proof.} By the local central limit theorem, $p_t(x) \leq C(1+t)^{-\frac{d}{2}}$ uniformly for $t>0$ and $x\in\Z^d$.
Hence
\be
\limsup_{t\to\infty} \frac{\sum_{x\in\Z} p_{\rho t}(x)\log p_t(x)}{\log t} \leq -\frac{d}{2}.
\ee
For a matching lower bound, we need lower bounds for $p_{t}(x)$ for all $x\in\Z^d$. Note that if $p^{(1)}_t(\cdot)$
denotes the transition probability kernel of a rate 1 simple random walk on $\Z$, then $p_t(x) = \Pi_{i=1}^d
p^{(1)}_{t/d}(x_i)$, and $\sum_{x\in\Z^d} p_{\rho t}(x) \log p_{t}(x) = d
\sum_{x\in\Z} p^{(1)}_{\rho t/d}(x) \log p^{(1)}_{t/d}(x)$. Hence it suffices to show
\be\label{128}
\liminf_{t\to\infty} \frac{\sum_{x\in\Z} p^{(1)}_{\rho t}(x) \log p^{(1)}_{t}(x)}{\log t} \geq -\frac{1}{2}.
\ee
For $0<\epsilon <\!\!<1 <\!\!< A <\infty$, we have the following estimates. There exist $C_1, C_2,C_3, T>0$ depending
on $\epsilon$ and $A$, such that
\begin{eqnarray}
\label{pest2}
p^{(1)}_t(x) &\geq& C_1 t^{-\frac{1}{2}} e^{-\frac{C_2x^2}{t}} \qquad\ \ \, \forall\ t\geq T,\ |x|\leq \epsilon t, \\
\label{pest3}
p^{(1)}_t(x) &\geq& e^{-C_3t} \qquad \qquad \qquad\, \forall\ t\geq T,\ \eps t< |x|< At, \\
\label{pest4}
p_t^{(1)}(x) &\geq& e^{-2 |x|\log |x|} \qquad \qquad \forall\ t\geq T,\ At \leq |x|.
\end{eqnarray}
To derive (\ref{128}) from (\ref{pest2})--(\ref{pest4}), we partition the sum $\sum_{x\in\Z}$ into
$\sum_{|x|\leq \eps t}$, $\sum_{\eps t<|x| < A t}$, and $\sum_{|x|\geq A t}$ with
$\epsilon <\!\!< \rho <\!\!< A$. By (\ref{pest2}),
\ben
\!\!\!\!\!\!\!\!\!
\sum_{|x|\leq \epsilon t}  p^{(1)}_{\rho t}(x) \log p^{(1)}_{t}(x)\!\!\!&\geq&\!\!\! \sum_{|x|\leq \epsilon t} p^{(1)}_{\rho t}(x)
\log\big(C_1 t^{-\frac{1}{2}} e^{-\frac{C_2 x^2}{t}}\big) \nonumber \\
\!\!\!&\geq&\!\!\! -\frac{\log t}{2} - |\log C_1| - \frac{C_2}{t} \sum_{x\in\Z} x^2 p^{(1)}_{\rho t}(x)
= -\frac{\log t}{2} - |\log C_1| - C_2\rho.  \label{est1}
\een
By (\ref{pest3}) and the Markov inequality,
\be\label{est2}
\sum_{\epsilon t <|x|< A t}  p^{(1)}_{\rho t}(x) \log p^{(1)}_{t}(x) \geq -C_3 t \sum_{|x|> \epsilon t}p^{(1)}_{\rho t}(x)
\geq -C_3 t \frac{\sum_{x\in\Z}x^2p^{(1)}_{\rho t}(x)}{\epsilon^2 t^2} = -\frac{C_3\rho}{\epsilon^2}.
\ee
And by (\ref{pest4}), for $t$ sufficiently large, we have
\be\label{est3}
\sum_{|x|\geq A t}\!\!  p^{(1)}_{\rho t}(x) \log p^{(1)}_{t}(x) \geq -2\!\!\! \sum_{|x|\geq A t}\!\!\! p^{(1)}_{\rho
  t}(x) |x|\log|x| \geq -2\!\!\! \sum_{|x|\geq A t}\!\!\! p^{(1)}_{\rho t}(x) \frac{|x|^2}{\frac{A t}{\log(A t)}}
\geq -\frac{2\rho}{A} \log (A t).
\ee
Combining (\ref{est1})--(\ref{est3}), we obtain the lower bound
\be
\liminf_{t\to\infty} \frac{\sum_{x\in\Z}p^{(1)}_{\rho t}(x) \log p^{(1)}_{t}(x)}{\log t} \geq -\frac{1}{2}-\frac{2\rho}{A}.
\ee
Since $A$ can be chosen arbitrarily large, (\ref{128}) follows.

We now verify (\ref{pest2})--(\ref{pest4}). Let $P_n(x)$ denote the probability that a discrete time simple random walk
starting from 0 visits $x$ at time $n$. Then for $x$ and $n$ having the same parity, by Stirling's formula,
\ben
P_n(x) = \frac{1}{2^n} \frac{n!}{\left(\frac{n+x}{2}\right)!\left(\frac{n-x}{2}\right)!} &=& \frac{(1+o(1))\sqrt{2\pi
    n}\left(\frac{n}{e}\right)^n}{2^n \sqrt{2\pi\left(\frac{n+x}{2}\right)}\left(\frac{n+x}{2e}\right)^{\frac{n+x}{2}}
\sqrt{2\pi\left(\frac{n-x}{2}\right)}\left(\frac{n-x}{2e}\right)^{\frac{n-x}{2}}} \nonumber\\
&=& (1+o(1))\sqrt{\frac{2n}{\pi(n^2-x^2)}}\ \
e^{-\left(\frac{n+x}{2}\right)\log\left(1+\frac{x}{n}\right)-\left(\frac{n-x}{2}\right)\log\left(1-\frac{x}{n}\right)}
\nonumber \\
&=& (1+o(1)) \sqrt{\frac{2n}{\pi(n^2-x^2)}}\ \  e^{-\frac{x^2}{2n}+ o\left(\frac{x^2}{n^2}\right)n}.
\een
Hence for $n$ sufficiently large and $|x|/n$ sufficiently small, we have
\be\label{139}
P_n(x) \geq C n^{-\frac{1}{2}} e^{-\frac{x^2}{n}}.
\ee
If $N_t$ denotes a Poisson random variable with mean $t$, then (\ref{pest2}) follows from (\ref{139}) and
the observation that  $N_t/t \to 1$ in probability with $|\P(N_t {\rm\ is\ odd}) - \P(N_t {\rm\ is\ even})| \to 0$ as $t\to\infty$.

For (\ref{pest3}), note that for $|x| <At$, by (\ref{139}),
\ben
p^{(1)}_t(x) \geq \sum_{At/\epsilon \leq n \leq 2At/\epsilon \atop n\equiv x {\rm\ mod\ 2}} \P(N_t = n) P_n(x)
&\geq& C \sqrt{\frac{\epsilon}{2At}} e^{-\frac{\epsilon x^2}{At}} \P\big(At/\epsilon \leq N_t \leq 2At/\epsilon,\
  N_t\equiv x {\rm\ mod\ 2}\big) \nonumber\\
&\geq& C \sqrt{\frac{\epsilon}{2At}} e^{-\epsilon At}  e^{-C't} \geq e^{-C_3 t},
\een
where we used the fact that $N_t/t$ satisfies a large deviation principle with a finite rate function on $[0,\infty)$.

For $|x|\geq At$, we can bound $p^{(1)}_t(x)$ from below by requiring that the random walk makes exactly $|x|$ jumps in
the time interval $[0,1]$ so that the random walk is at $x$ at time 1, and at time $t$ the random walk returns to $x$.
Thus, by the local central limit theorem, for $t$ large,
\be
p^{(1)}_t(x) \geq \frac{1}{e |x|!} 2^{-|x|} \frac{C}{t} =
(1+o(1)) \frac{e^{-1 +|x| - |x|\log |x|}}{\sqrt{2\pi |x|}} 2^{-|x|} \frac{C}{t}.
\ee
It is then clear that (\ref{pest4}) holds. \qed
\medskip

\noindent
{\bf Remark.} We point out that, for general mean zero finite variance random walks, the estimates
(\ref{pest2})--(\ref{pest4}) can still be established by adapting the proof here and decomposing the
random walk transition kernel to extract a simple random walk part.
\medskip

\noindent
{\bf Remark.} The analogue of Lemma \ref{L:srwest} also holds for discrete time simple random walks. The proof
is similar and omitted.

\bigskip

\noindent
{\bf Lower bound for $\E^Y_0\left[\log \E^X_0\left[e^{\beta L_t(X,Y)} \big| X_t=Y_t\right] \right]$ for $d=1$:}

\noindent
By Jensen's inequality,
\be\nonumber
\E^Y_0\left[\log \E^X_0\big[e^{\beta L_t(X,Y)} \big | X_t=Y_t\big] \right]
\geq
\E^Y_0\left[\E^X_0\big[\beta L_t(X,Y)\big| X_t=Y_t\big]\right]
\!=\! \beta\!\! \int_0^t \E^Y_0\left[\frac{p_{s}(Y_s)p_{(t-s)}(Y_t-Y_s)}{p_{t}(Y_t)}\right] ds. \label{144}
\ee
By Donsker's invariance principle, there exists $\alpha>0$ s.t.\ $\P^Y_0(\sup_{s\in [0,t]}|Y_s|\leq\sqrt{t})\geq\alpha$
for all $t>0$. On the other hand, if $\sup_{s\in [0,t]}|Y_s|\leq\sqrt{t}$, then by the local central limit theorem,
$p_s(Y_s)\wedge p_{t-s}(Y_t-Y_s) \geq C/\sqrt{t}$ for all $s\in [t/3, 2t/3]$ for some $C$ independent of $Y$ and
$t>1$, while $p_t(Y_t) \leq C'/\sqrt{t}$. Therefore
\ben
\E^Y_0\left[\log \E^X_0\big[e^{\beta L_t(X,Y)} \big | X_t=Y_t\big] \right] \geq
 \alpha \beta \int_{t/3}^{2t/3} \frac{\frac{C}{\sqrt t}\frac{C}{\sqrt t}}{\frac{C'}{\sqrt t}}ds = C'\sqrt{t}
\een
for some $C'>0$ independent of $t$. In view of (\ref{125}) and Lemma \ref{L:srwest}, this proves that
$\E^Y_0[\log Z^{\beta,{\rm pin}}_{t,Y}]>0$ for $t$ large, and hence $\beta_{\rm c}=0$ for $d=1$.
\bigskip

\noindent
{\bf Lower bound for $\E^Y_0\left[\log \E^X_0\big[e^{\beta L_t(X,Y)}\,\big| X_t=Y_t\big] \right]$ for $d=2$:}

\noindent
Since in $d=2$, $L_t(X,Y)$ is typically of order $\log t$, the argument above for $d=1$ fails for $d=2$.
Instead, we apply an a.s.\ limit theorem for $L_t(X,Y)/\log t$ conditioned on $Y$. More precisely, by
Theorem 1.2 of G\"artner and Sun \cite{GS07}, a.s.\ w.r.t.\ $Y$, $L_t(X,Y)/\log t$ conditioned on $Y$
converges in distribution to an exponential random variable with mean $1/\pi(1+\rho)$. We only need to bypass
the conditioning on $X_t=Y_t$.

Let $\mu_{t/\log t}$ denote the law of $(X_s)_{0\leq s\leq t/\log t}$, and let $\mu^{(t,y)}_{t/\log t}$ denote
the law of $(X_s)_{0\leq s\leq t/\log t}$ conditioned on $X_t=y$. Then $\mu_{t/\log t}$ and $\mu^{(t,y)}_{t/\log t}$ are
equivalent with density
\be
\frac{d\,\mu^{(t,y)}_{t/\log t}}{d\,\mu_{t/\log t}}\big(X\big) =
\frac{p_{t-t/\log t}(y-X_{t/\log t})}{p_{t}(y)} =
\frac{t}{t-t/\log t}
\frac{e^{-\frac{\Vert y-X_{t/\log t}\Vert^2} {t-t/\log t} } +o(1)}{e^{ -\frac{\Vert y\Vert^2}{t} }+o(1)},
\ee
where we applied the local central limit theorem. Since $\Vert X_{t/\log t}\Vert/\sqrt{t}\to 0$ in probability
as $t\to\infty$, it is clear that in total variational distance,
\be
\sup_{\Vert y\Vert \leq \sqrt{t}} \big\Vert \mu^{(t,y)}_{t/\log t}- \mu_{t/\log t}\big\Vert_{\rm TV}
\underset{t\to\infty}{\longrightarrow} 0.
\ee
We can thus remove the conditioning at the cost of reducing the time interval from $t$ to $t/\log t$.

Fix $A>0$. Let
\be
G^A_{t/\log t} = \big\{Y : \mu_{t/\log t}\big(L_{t/\log t}(X,Y) \geq A\log t\big) \geq e^{-\alpha A}\big\}.
\ee
By Theorem 1.2 of \cite{GS07}, if we choose $\alpha >\pi(1+\rho)$, then $\P^Y_0(G^A_{t/\log t}) \to 1$ as
$t\to \infty$. We now write
\begin{eqnarray}
&& \E^Y_0\left[\log \E^X_0\left[e^{\beta L_t(X,Y)} \Big| X_t=Y_t\right] \right] \nonumber\\
&\geq& \E^Y_0\left[1_{\{\Vert Y_t\Vert \leq \sqrt{t},\ Y\in G^A_{t/\log t}\}} \log \E^X_0\left[e^{\beta L_{t/\log t}(X,Y)}\Big|
    X_t=Y_t\right]\right] \nonumber \\
&\geq& \E^Y_0\left[1_{\{\Vert Y_t\Vert \leq \sqrt{t},\ Y\in G^A_{t/\log t}\}}
\Big(\beta A\log t + \log \mu^{(t,Y_t)}_{t/\log t}\big(L_{t/\log t}(X,Y)\geq
  A\log t\big)\Big)\right] \nonumber \\
&\geq& \beta A\, \P^Y_0(\Vert Y_t\Vert \leq\sqrt{t},\, Y\in G^A_{t/\log t}) \log t \nonumber \\
&& \qquad \qquad \qquad +\ \E^Y_0\left[1_{\{\Vert Y_t\Vert \leq
   \sqrt{t},\, Y\in G^A_{t/\log t}\}}\log \Big(\mu_{t/\log t}\big(L_{t/\log t}(X,Y)\geq  A\log t\big)+o(1)\Big)\right].
\nonumber\\
&\geq& (C-o(1))(\beta A \log t + \log(e^{-\alpha A} +o(1)),
\end{eqnarray}
where $C= \inf_{t>0}\P^Y_0(\Vert Y_t\Vert\leq \sqrt{t})$ is positive and independent of $A$. Since $A$ can be chosen
arbitrarily large, in view of (\ref{125}) and Lemma \ref{L:srwest}, this proves that $\E^Y_0[\log Z^{\beta,{\rm pin}}_{t,Y}]>0$
for $t$ large, and hence $\beta_{\rm c}=0$ for $d=2$.
\qed

\section{Gap between critical points: discrete time}
\subsection{Proof of Theorem \ref{T:cptst} in discrete time: $d\geq 5$}
Our proof is based on adaptations of the fractional moment method used recently by Derrida, Giacomin, Lacoin
and Toninelli \cite{DGLT07} to show the non-coincidence of annealed and quenched critical points for the pinning
model in the disorder-relevant regime. Two ingredients are needed for the adaptation. First, a suitable representation
for the partition function $\hat Z^\beta_{N,Y}$ and its constrained counterpart $\hat Z^{\beta, {\rm pin}}_{N,Y}$ in a
similar form as in (\ref{hatzsplit2}), except with a Gibbs weight factor $w(\cdot)$ that has a simpler dependence on the
disorder $(\Delta_i)_{i\in\N}=(Y_{i+1}-Y_i)_{i\in\N}$ than in (\ref{rwf}). Second, a suitable change of measure for the
disorder $Y$ when estimating fractional moments $\E^Y_0[(\hat Z^{\beta, {\rm pin}}_{N,Y})^\gamma]$ for $N$ on the order of
the correlation length of the annealed model.

We split the proof into three parts: representation for $\hat Z^\beta_{N,Y}$ and $\hat Z^{\beta, {\rm pin}}_{N,Y}$;
fractional moment method; change of measure. To simplify notation, $C, C_1, C'$, etc, will denote generic constants
whose precise values may change from place to place.
\bigskip

\noindent
{\bf Representation for $\hat Z^\beta_{N,Y}$ and $\hat Z^{\beta, {\rm pin}}_{N,Y}$.} The representation we now
derive was already used in \cite{BGdH08}. It is based on binomial expansion for $(1+e^\beta-1)^{L_N(X,Y)}$. Let
$p^X_n(\cdot)$, resp.\ $p^{X-Y}_n(\cdot)$, be the $n$-step transition probability kernel of $X$, resp.\ $X-Y$. Let
$G^{X-Y} = \sum_{n=1}^\infty p^{X-Y}_n(0)$, $K(n) = p^{X-Y}_n(0)/G^{X-Y}$, $z'= e^\beta-1$, $z=z'G^{X-Y}$, and
$\check Z^z_{N,Y} = \hat Z^\beta_{N,Y}$. Then
\begin{eqnarray}
\check Z^z_{N,Y}
&=& \E^X_0\left[(1+z')^{L_N(X,Y)}\right]
= \E^X_0\Big[1+\sum_{m=1}^{N}\sum_{\sigma_0=0<\sigma_1<\cdots <\sigma_{m}\leq N} (z')^m\prod_{i=1}^{m} 1_{\{X_{\sigma_i}=Y_{\sigma_i}\}}\Big] \nonumber \\
&=& 1+\sum_{m=1}^{N} \sum_{\sigma_0=0<\sigma_1<\cdots <\sigma_m\leq N}\!\!\!\!\!\!
(z')^m \prod_{i=1}^m p^X_{\sigma_i-\sigma_{i-1}}(Y_{\sigma_i}-Y_{\sigma_{i-1}}) \nonumber \\
&=& 1+ \sum_{m=1}^{N} \sum_{\sigma_0=0<\sigma_1<\cdots <\sigma_m\leq N}
\prod_{i=1}^m K(\sigma_i-\sigma_{i-1}) w\big(z, \sigma_{i}-\sigma_{i-1}, Y_{\sigma_{i}}-Y_{\sigma_{i-1}}\big), \label{diszfree}
\end{eqnarray}
where
\be\label{wdiscrete}
w(z,\sigma_i-\sigma_{i-1}, Y_{\sigma_i}-Y_{\sigma_{i-1}})=
z p^X_{\sigma_i-\sigma_{i-1}}(Y_{\sigma_i}-Y_{\sigma_{i-1}})/p^{X-Y}_{\sigma_i-\sigma_{i-1}}(0).
\ee
If we denote $\check Z^{z,\rm pin}_{N,Y}=\frac{z'}{1+z'}\hat Z^{\beta, \rm pin}_{N,Y}$, then similarly,
\begin{eqnarray}
\check Z^{z, {\rm pin}}_{N,Y}
&=& \E^X_0\left[(1+z')^{L_{N-1}(X,Y)}z' 1_{\{X_N=Y_N\}}\right] \nonumber \\
&=& \sum_{m=1}^{N} \sum_{\sigma_0=0<\sigma_1<\cdots <\sigma_m=N}
\prod_{i=1}^m K(\sigma_i-\sigma_{i-1}) w\big(z, \sigma_{i}-\sigma_{i-1}, Y_{\sigma_{i}}-Y_{\sigma_{i-1}}\big). \label{disz}
\end{eqnarray}
Note that (\ref{disz}) casts $\check Z^{z, {\rm pin}}_{N,Y}$ in the same form as (\ref{hatzsplit2}), except now
$K(n)$ equals $p^{X-Y}_n(0)/G^{X-Y}$ instead of $\P^{X-Y}(\tau_0=n)$.
This mapping from one underlying renewal process to another defined in terms of the Green function
decomposition of the original renewal process applies to any pinning model with an underlying transient renewal distribution.
Of course the disorder also changes and the terms in (\ref{disz}) may not be positive in general.
This is not the case here, and the key point for us is that the weight factor $w$ now has a much simpler dependence on
the disorder $(\Delta_j)_{\sigma_{i-1}<j\leq\sigma_i}$ (i.e.\ only on $\sigma_i-\sigma_{i-1}$ and
$\sum_{j=\sigma_{i-1}+1}^{\sigma_i}\Delta_j$) than in (\ref{rwf}). We note that if $\Tilde K(n) \sim \frac{c}{n^{1+\alpha}}$
for some $\alpha>0$ is the first return time distribution of a transient renewal process, then the corresponding return probability
at time $n$ satisfies $p(n)\sim \frac{c'}{n^{1+\alpha}}$. See \cite[Theorem A.4]{G07}.

Because $K$ is the return time
distribution of a recurrent renewal process $\sigma$ on $\N_0$, and
$\E^Y_0[w(z, \sigma_i-\sigma_{i-1},Y_{\sigma_i}-Y_{\sigma_{i-1}})]=z$, the critical point for the annealed model
associated with $\check Z^{z, \rm pin}_{N,Y}$  is $z^{\rm ann}_{\rm c}=1$, or equivalently,
$1=z^{\rm ann}_{\rm c}=(e^{\hat \beta^{\rm ann}_{\rm c}}-1)G^{X-Y}$ so that
\be\label{abcdis}
\hat \beta^{\rm ann}_{\rm c} = \log\left(1+\frac{1}{G^{X-Y}}\right).
\ee

\noindent
{\bf Fractional moment method.}
We now recall the fractional moment method used by Derrida et al in \cite{DGLT07}. Due to the common framework
between pinning models and the random walk pinning model as pointed out in Section 1.3, the basic strategy
carries over without change. The only model dependent part of the argument lies in estimating
$\E^Y_0[(\check Z^{z,{\rm pin}}_{N,Y})^\gamma]$, $\gamma\in(0,1)$, for $N$ on the order of the correlation length of the annealed
model, where a change of measure argument for the disorder needs to be adapted.

In terms of the new variables $z=(e^\beta-1)G^{X-Y}$ and $\check Z^z_{N,Y}$, Theorem \ref{T:cptst} reduces to
showing that for some $z>z^{\rm ann}_{\rm c}=1$, $\sup_{N\in\N_0} \check Z^z_{N,Y}<\infty$ a.s.\ w.r.t.\ $Y$.
Since for $z>1$, $\check Z^z_{N,Y}$ is a.s.\ increasing in $N$, it suffices to show that for some $z>1$ and
$\gamma\in (0,1)$,
\be\label{155}
\sup_{N\in\N_0} \E^Y_0\big[\big(\check Z^{z}_{N,Y}\big)^\gamma\big] < \infty.
\ee
The basic idea is to suitably group terms in the expansion for $\check Z^{z}_{N,Y}$ in (\ref{diszfree})
and then apply the fractional moment inequality
\be\label{fracineq}
\Big(\sum_{i=1}^n |a_i|\Big)^\gamma \leq \sum_{i=1}^n |a_i|^\gamma, \qquad \gamma \in (0,1).
\ee
However, the effectiveness of (\ref{fracineq}) depends crucially on how $\check Z^{z}_{N,Y}$ is
decomposed. In \cite{DGLT07}, Derrida et al studied analogues  of the constrained partition
function $\check Z^{z, \rm pin}_{N,Y}$, and their clever choice is to group terms in (\ref{disz}) according
to the starting and the ending position of the gap in the renewal sequence $\sigma$ straddling a fixed position
$L\in\N$. Namely,
\be\nonumber
\check Z^{z, {\rm pin}}_{N,Y} = \sum_{i=0}^{L-1}\sum_{j=0}^{N-L}
\check Z^{z, {\rm pin}}_{i,Y} K(N-j-i)w(z, N-j-i,Y_{N-j}-Y_i)
\check Z^{z, {\rm pin}}_{j, \theta_{N-j}Y},
\ee
where $\theta_nY = (Y_{n+i}-Y_n)_{i\in\N_0}$ denotes a shift in $Y$. For $\check Z^z_{N,Y}$, we can perform
a similar grouping of terms in (\ref{diszfree}) and get
\be\label{checkZsplit}
\check Z^z_{N,Y} = \check Z^z_{L-1, Y} + \sum_{i=0}^{L-1} \sum_{j=0}^{N-L} \check Z^{z,\rm pin}_{i,Y}
K(N-j-i)w(z,N-j-i,Y_{N-j}-Y_i) \check Z^z_{j, \theta_{N-j}Y}.
\ee
Fix $\gamma \in (0,1)$. Denote $\check A^z_N = \E^Y_0\big[(\check Z^z_{N,Y})^\gamma\big]$ and
$\check A^{z,\rm pin}_N=\E^Y_0\big[(\check Z^{z,{\rm pin}}_{N,Y})^\gamma\big]$. Since
$$
K(N-j-i)w(z, N-j-i,Y_{N-j}-Y_i)=\frac{z p^X_{N-j-i}(Y_{N-j}-Y_i)}{G^{X-Y}}  \leq C(N-j-i)^{-\frac{d}{2}}
$$
for some $C>0$ independent of $i$, $j$, $N$, $Y$ and $z \in [1,2]$ by the local central limit theorem, applying
(\ref{fracineq}) to (\ref{checkZsplit}) and taking expectation w.r.t.\ $Y$ gives
\be\label{453}
\check A^z_N \leq \check A^z_{L-1} +
C \sum_{i=0}^{L-1} \check A^{z,\rm pin}_i \sum_{j=0}^{N-L} (N-j-i)^{-\frac{d\gamma}{2}} \check A^{z}_j
\leq \check A^z_{L-1} +
C\left(\sum_{i=0}^{L-1}\frac{\check A^{z,\rm pin}_i}{(L-i)^{\frac{d\gamma}{2}-1}}\right) \max_{0\leq j\leq
  N-L}\check A^{z}_j.
\ee
If for some choice of $z>1$ and $L\in\N$,
\be\label{rhodis}
\check\varrho = C\left(\sum_{i=0}^{L-1}\frac{\check A^{z,\rm pin}_i}{(L-i)^{\frac{d\gamma}{2}-1}}\right) <1,
\ee
then iterating (\ref{453}) clearly implies that $\check A^{z}_N$ is uniformly bounded in $N$, and hence (\ref{155}).

By Jensen's inequality, $\check A^{z,\rm pin}_N \leq \E^Y_0[\check Z^{z,\rm pin}_{N,Y}]^\gamma$. It is clear from
(\ref{disz}) and (\ref{wdiscrete}) that $\E^Y_0[\check Z^{z,\rm pin}_{N,Y}]$ is the partition function of a
homogeneous pinning model with critical point $z^{\rm ann}_{\rm c}=1$. Hence
$\check F_{\rm ann}(z)= \displaystyle\lim_{N\to\infty}N^{-1}\log \E^Y_0[\check Z^{z,{\rm pin}}_{N,Y}]$ exists, and
$\check F_{\rm ann}(z)=\hat F_{\rm ann}(\beta)$ with $z=(e^\beta-1)G^{X-Y}$. Since $d\geq 5$, $K(\cdot)$ has finite first
moment, and hence by Theorem 2.1 of \cite{G07}, $\check F_{\rm ann}(z) \sim C(z-1)$ for some $C>0$ as $z\downarrow 1$.
Since $(\E^Y_0[\check Z^{z, \rm pin}_{n,Y}])_{n\in\N}$ is super-multiplicative,
$\E^Y_0[\check Z^{z, \rm pin}_{N,Y}] \leq e^{N\check F_{\rm ann}(z)}\leq e^{CN(z-1)}$ for all $N\in\N$. So if we choose
\be\label{discreteLz}
L=L(z)=\frac{1}{z-1},
\ee
where we abused notation and assumed $L$ to be an integer for simplicity, then
$\sup_{1\leq i\leq L}\check A^{z,\rm pin}_i\leq C$ for some $C>0$ independent of $z$. Therefore
\be\label{456}
\check\varrho \leq \sum_{i=0}^{L-R} \frac{C}{(L-i)^{\frac{d\gamma}{2}-1}} + \sum_{i=L-R+1}^{L-1}\frac{C \check A^{z,\rm
    pin}_i}{(L-i)^{\frac{d\gamma}{2}-1}} \leq C R^{2-\frac{d\gamma}{2}} + C\max_{L-R\leq i\leq L}\check A^{z,\rm pin}_i.
\ee
For $d\geq 5$, we can choose $\gamma<1$ close to $1$ such that the first term on the RHS of (\ref{456}) can be made
arbitrarily small (uniformly in $z$) by choosing $R$ large. To show $\check\varrho< 1$ for some $z>1$, it then suffices
to show that
\be\label{ALbound}
\lim_{z\downarrow 1} \max_{L-R\leq N\leq L} \check A^{z,\rm pin}_N=0,
\ee
where $R\in\N$ is large and fixed, and $L=\frac{1}{z-1}$. This summarizes the model independent part of the fractional
moment method as used in \cite{DGLT07}.
\bigskip

\noindent
{\bf Change of measure.}
The basic idea in \cite{DGLT07} to prove (\ref{ALbound}) is to apply a change of measure to the disorder so that
the cost of changing the measure is small, yet under the new disorder, the annealed partition function for a
system of size $L$ is small. For the pinning model, the choice of changing the measure in \cite{DGLT07} is to
make the disorder more repulsive, i.e., tilt the measure of $\omega_i$ in (\ref{115}) by a factor
$e^{-\lambda\omega_i}$ for some $\lambda>0$. In our setting, it turns out that for the continuous time model,
the appropriate change of measure is to increase the jump rate of the random walk $Y$. For the discrete time
model, the analogue is to increase the variance of the random walk increment each step without changing the
support of the random walk transition kernel. However, among nearest-neighbor random walks on $\Z^d$, the
variance of simple random walk is already maximal. To overcome this difficulty,  we change measure for $Y$ two
steps at a time. More precisely, for $h\in (0, \frac{1}{2d})$, let $(Y^h_n)_{n\in\N_0}$ be a process
on $\Z^d$ with $Y_0=0$ and transition probabilities
\be\label{Yhtransit}
\P\big(Y^h_{n+1}-Y^h_n =e_i \big| (Y^h_k)_{0\leq k\leq n}\big) =
\left\{ \begin{aligned}
\frac{1}{2d}\quad   & \qquad \mbox{ if $n$ is even, or $n$ is odd and $e_i\neq \pm (Y^h_n- Y^h_{n-1})$}, \\
\frac{1+h}{2d} & \qquad \mbox{ if $n$ is odd, and $e_i= Y^h_n-Y^h_{n-1}$}, \\
\frac{1-h}{2d} & \qquad \mbox{ if $n$ is odd, and $e_i= -(Y^h_n-Y^h_{n-1})$},
\end{aligned}
\right.
\ee
for each of the $2d$ unit vectors $e_i\in\Z^d$. Note that $\P(Y^h_2=2e_i) = \P(Y_2=2e_i)+\frac{h}{4d^2}$ for
each unit vector $e_i\in\Z^d$, $\P(Y^h_2=0)= \P(Y_2=0)-\frac{h}{2d}$, and $\P(Y^h_2=x)=\P(Y_2=x)$ for
all other $x\in\Z^d$. Thus $Y^h_2$ has larger variances than $Y_2$. Clearly up to any time $N\in\N$, the
distribution of $Y$ and $Y^h$ are equivalent. Let $f(N,Y)$ denote the Radon-Nikodym derivative of the law
of $(Y^h_i)_{0\leq i\leq N}$ w.r.t.\ $(Y_i)_{0\leq i\leq N}$. Then
\begin{eqnarray}
\nonumber
\check A^{z,\rm pin}_N = \E^{Y^h}_0\big[f(N,Y^h)^{-1}\big(\check Z^{z,\rm pin}_{N,Y^h}\big)^\gamma \big]
&\leq& \E^{Y^h}_0\big[f(N,Y^h)^{-\frac{1}{1-\gamma}}\big]^{1-\gamma} \E^{Y^h}_0\big[\check Z^{z,\rm pin}_{N,Y^h}\big]^\gamma\\
&=& \E^Y_0\big[f(N,Y)^{-\frac{\gamma}{1-\gamma}}\big]^{1-\gamma}\,
\E^{Y^h}_0\big[\check Z^{z,\rm pin}_{N,Y^h}\big]^\gamma. \label{185}
\end{eqnarray}
Since $(Y_{2n+1}-Y_{2n},Y_{2n+2}-Y_{2n})_{n\in\N_0}$ are i.i.d.\ and the distribution of $Y^h_{2n+1}-Y^h_{2n}$ conditioned
on $Y^h_{2n}$ is the same as a simple random walk, we have
\begin{eqnarray}
\E^Y_0\big[f(N,Y)^{-\frac{\gamma}{1-\gamma}}\big] =
\E^Y_0\big[f(2,Y)^{-\frac{\gamma}{1-\gamma}}\big]^{\left\lfloor\frac{N}{2}\right\rfloor}
= \Big(1-\frac{1}{d}+\frac{(1+h)^{-\frac{\gamma}{1-\gamma}}}{2d}
+\frac{(1-h)^{-\frac{\gamma}{1-\gamma}}}{2d}\Big)^{\left\lfloor\frac{N}{2}\right\rfloor}
\leq e^{\frac{\gamma h^2 N}{2d(1-\gamma)^2}} \nonumber
\end{eqnarray}
for $h$ sufficiently small. Therefore if we choose $h=\frac{1}{\sqrt L}$, then the first factor in (\ref{185}) is
uniformly bounded for $L-R\leq N\leq L$, and to prove (\ref{ALbound}), it only remains to estimate
$\E^{Y^h}_0[\check Z^{z, \rm pin}_{N,Y^h}]$ for $h=\frac{1}{\sqrt L}=\sqrt{z-1}$.

By (\ref{disz}), we have
\be\label{460}
\E^{Y^h}_0\big[\check Z^{z,\rm pin}_{N,Y^h}] = \sum_{m=1}^{N} \Big(\frac{z}{G^{X-Y}}\Big)^{m}\!\!\!\!\!\!
\sum_{\sigma_0=0<\sigma_1<\cdots<\sigma_m=N}\!\!\!\!
\E^{Y^h}_0\Big[\prod_{i=1}^m p^X_{\sigma_i-\sigma_{i-1}}(Y^h_{\sigma_i}-Y^h_{\sigma_{i-1}})\Big].
\ee
Note that when $\sigma_{i-1}$ is even, by the properties of $Y^h$, we have
$$
\E^{Y^h}_0\big[p^X_{\sigma_i-\sigma_{i-1}}(Y^h_{\sigma_i}-Y^h_{\sigma_{i-1}})\big| (Y^h_j)_{0\leq j\leq \sigma_{i-1}}\big]
=\E^{Y^h}_0\big[p^X_{\sigma_i-\sigma_{i-1}}(Y^h_{\sigma_i-\sigma_{i-1}})\big].
$$
Similarly when $\sigma_{i-1}$ is odd, by symmetry and translation invariance, we have
$$
\E^{Y^h}_0\big[p^X_{\sigma_i-\sigma_{i-1}}(Y^h_{\sigma_i}-Y^h_{\sigma_{i-1}})\big| (Y^h_j)_{0\leq j\leq \sigma_{i-1}}\big]
= \E^{Y^h}_0\big[p^X_{\sigma_i-\sigma_{i-1}}(Y^h_{\sigma_i-\sigma_{i-1}+1}-Y^h_1)\big| Y^h_1=e_1\big],
$$
which is a constant independent of $(Y^h_j)_{0\leq j\leq \sigma_{i-1}}$. Thus in (\ref{460}), we can successively
condition w.r.t.\ $(Y^h_j)_{0\leq j\leq \sigma_{n}}$, $(Y^h_j)_{0\leq j\leq \sigma_{n-1}}$, $\ldots$,
$(Y^h_j)_{0\leq j\leq \sigma_{1}}$. To write the result in a more compact form, let us denote
\be
\begin{aligned}
K_{h,\rm even}(n) = \frac{\E^{Y^h}_0[p^X_n(Y^h_n)]}{G_{h,\rm even}} & \qquad \mbox{where} \quad
G_{h,\rm even} = \sum_{n=1}^\infty \E^{Y^h}_0[p^X_n(Y^h_n)], \\
K_{h,\rm odd}(n) = \frac{\E^{Y^h}_0[p^X_n(Y^h_{n+1}-Y^h_1) | Y^h_1=e_1]}{G_{h,\rm odd}} & \qquad \mbox{where} \quad
G_{h, \rm odd} = \sum_{n=1}^\infty \E^{Y^h}_0[p^X_{n}(Y^h_{n+1}-Y^h_1) | Y^h_1=e_1].
\end{aligned} \nonumber
\ee
Let $K_h(i,j)=K_{h,\rm even}(j-i)$ when $i$ is even, and $K_h(i,j)=K_{h, \rm odd}(j-i)$ when $i$ is odd.
Let $\iota=\{0, \iota_1, \iota_2, \cdots\}$ be a renewal process on $\N_0$ with parity-dependent inter-arrival law
$K_h(\cdot,\cdot)$, and denote expectation w.r.t.\ $\iota$ by $\E^{K_h}[\cdot]$. Then (\ref{460}) reduces to
$$
\E^{Y^h}_0\big[\check Z^{z,\rm pin}_{N,Y^h}\big] =
\E^{K_h}\Big[\Big(\frac{z}{G^{X-Y}}\Big)^{|\iota\cap[1,N]|} G_{h,\rm even}^{|\iota_{\rm e}\cap[1,N]|} G_{h,\rm odd}^{|\iota_{\rm o}\cap[1,N]|} 1_{\{N\in\iota\}}
 \Big] \leq \E^{K_h}\Big[\Big(\frac{z(G_{h,\rm even}\vee G_{h,\rm odd})}{G^{X-Y}}\Big)^{|\iota\cap [1,N]|}\Big],
$$
where $\iota_{\rm e}$ and $\iota_{\rm o}$ denote respectively the even and odd subsets of $\iota$. In $d\geq 5$,
by the local central limit theorem, it is easy to
see that there exists an inter-arrival probability distribution $K_*(\cdot)$ on $\N$ with finite first moment, such that
$K_*$ stochastically dominates both $K_{h,\rm even}(\cdot)$ and $K_{h, \rm odd}(\cdot)$ for $h$ sufficiently small,
i.e., $\sum_{i\geq n} K_*(i) \geq \sum_{i\geq n} K_{h,\rm even}(i)$ and
$\sum_{i\geq n} K_*(i) \geq \sum_{i\geq n} K_{h,\rm odd}(i)$ for all $n\in\N$ and $h\in [0,\frac{1}{2}]$.
Recall our choice $h=\frac{1}{\sqrt L}=\sqrt{z-1}$. We will show that
\be\label{187}
\frac{z\big(G_{h,\rm even}\vee G_{h,\rm odd}\big)}{G^{X-Y}}
= 1- c\sqrt{z-1}+o(\sqrt{z-1})
\ee
for some $c>0$. Then for all $z>1$ sufficiently close to 1,
\be\label{462}
\E^{Y^h}_0\big[\check Z^{z,\rm pin}_{N,Y^h}\big] \leq \E^{K_*}\big[\big(1- c\sqrt{z-1}+o(\sqrt{z-1})\big)^{|\iota^*\cap [1,N]|}\big],
\ee
where $\iota^*$ is a renewal process with inter-arrival law $K_*$ and is independent of $z$. By the law of large
numbers, a.s.\ w.r.t.\ $\iota^*$,
$$
\displaystyle \lim_{n\to\infty} N^{-1} |\iota^*\cap [1,N]| = \frac{1}{\sum_{i\in\N} iK_*(i)} >0,
$$
and hence
$$
\lim_{z\downarrow 1} \max_{(z-1)^{-1}-R\leq N\leq (z-1)^{-1}} \big(1- c\sqrt{z-1}+o(\sqrt{z-1})\big)^{|\iota^*\cap [1,N]|}=0.
$$
Thus
\be\label{463}
\lim_{z\downarrow 1} \max_{L-R\leq N\leq L} \E^{Y^h}_0\big[\check Z^{z,\rm pin}_{N,Y^h}\big]=0, \qquad L=\frac{1}{z-1},
\ h=\sqrt{z-1},
\ee
which together with (\ref{185}) implies (\ref{ALbound}).

It only remains to verify (\ref{187}). For $k=(k_1, \cdots, k_d)\in \R^d$, we have
\be
\begin{aligned}
\phi(k) &:= \E^X_0[e^{ik\cdot X_1}] = \frac{1}{d}\sum_{i=1}^d \cos k_i, \\
\psi(k) &:= \E^{Y^h}_0[e^{ik\cdot Y^h_2}] = \phi(k)^2 - \frac{h}{d^2}\sum_{i=1}^d \sin^2k_i, \\
\varphi(k) &:= \E^{Y^h}_0[e^{ik\cdot(Y^h_2-Y^h_1)} \,|\, Y^h_1=e_1] = \phi(k) + i\frac{h}{d}\sin k_1.
\end{aligned}
\ee
Since $X$ and $Y^h$ are independent, $(Y^h_{2n}-Y^h_{2n-2})_{n\in\N}$ are i.i.d., $Y^h_{2n+1}-Y^h_{2n}$ is independent
of $(Y^h_j)_{0\leq j\leq 2n}$ and is distributed as $X_1$, while conditioned on $Y^h_1=e_1$,
$Y^h_2-Y^h_1$ is independent of $(Y^h_j-Y^h_2)_{j\geq 2}$, we obtain by Fourier inversion
\begin{eqnarray}
G^{X-Y}&=& \frac{1}{(2\pi)^d} \int_{[-\pi,\pi]^d} \Big(\phi(k)^2+\phi(k)^4+\cdots\Big)dk
= \frac{1}{(2\pi)^d} \int_{[-\pi,\pi]^d} \frac{\phi(k)^2}{1-\phi(k)^2}dk,                 \label{189} \\
G_{h,\rm even} &=& \frac{1}{(2\pi)^d} \int_{[-\pi,\pi]^d} \Big(\phi(k)^2+
\phi(k)^2\psi(k)+\phi(k)^4\psi(k)+\phi(k)^4\psi(k)^2+\cdots \Big)dk                          \nonumber\\
&=& \frac{1}{(2\pi)^d} \int_{[-\pi,\pi]^d} \frac{\phi(k)^2(1+\psi(k))}{1-\phi(k)^2\psi(k)} dk, \label{190} \\
G_{h,\rm odd} &=& \frac{1}{(2\pi)^d} \int_{[-\pi,\pi]^d} \Big(\varphi(k)\phi(k)+
\varphi(k)\phi(k)^3+ \varphi(k)\phi(k)^3\psi(k) +\varphi(k)\phi(k)^5\psi(k) +\cdots \Big)dk \nonumber\\
&=& \frac{1}{(2\pi)^d} \int_{[-\pi,\pi]^d} \frac{\varphi(k)\phi(k)(1+\phi(k)^2)}{1-\phi(k)^2\psi(k)} dk
= \frac{1}{(2\pi)^d} \int_{[-\pi,\pi]^d} \frac{\phi(k)^2(1+\phi(k)^2)}{1-\phi(k)^2\psi(k)} dk, \label{191}
\end{eqnarray}
where in (\ref{191}) we have used the formula for $\varphi(k)$ and the fact that $\phi(k)$ and $\psi(k)$
are even functions while $\sin k_1$ is odd. Since $\psi(k) < \phi(k)^2$ and $\phi(k), \psi(k)\in[-1,1]$,
we have $G_{h,\rm even} < G_{h,\rm odd}$, while
\begin{eqnarray}
G^{X-Y}-G_{h,\rm odd} &=&
\frac{1}{(2\pi)^d} \int_{[-\pi,\pi]^d} \Big(\frac{\phi(k)^2}{1-\phi(k)^2}-
\frac{\phi(k)^2(1+\phi(k)^2)}{1-\phi(k)^2\psi(k)}\Big) dk, \nonumber\\
&=& \frac{h}{(2\pi)^d\, d^2} \int_{[-\pi,\pi]^d} \frac{\phi(k)^4\sum_{i=1}^d\sin^2k_i}{(1-\phi(k)^2)(1-\phi(k)^2\psi(k))}dk,
\end{eqnarray}
which implies (\ref{187}) since $h=\sqrt{z-1}$.
\qed
\medskip

\noindent
{\bf Remark.} Equation (\ref{187}) reveals the close resemblance between the random walk pinning model and
the pinning model (compare (\ref{462}) here with (4.12) in \cite{DGLT07})). In both cases, after changing
the measure, we end up comparing with a homogeneous pinning model of size $N$ with weight factor $e^{-c/\sqrt{N}}$
for each renewal return. The factor $c/\sqrt{N}$ partly explains why $\alpha=1/2$, resp.\ $d=3$, is the
critical case for the pinning, resp.\ random walk pinning model.
\medskip

\noindent
{\bf Remark.} For general random walks, we can try to change measure for $Y$ one-step at a time. More precisely, let
$S=\{y\in\Z^d : p^Y_1(y)>0\}$. Then for any $A, B\subset S$ and for any transition probability
kernels $p^A_1(\cdot)$ and $p^B_1(\cdot)$ with support resp.\ $A$ and $B$, and for $h\in \R$ sufficiently close
to $0$, we can change measure for $Y$ by replacing $p^Y_1(\cdot)$ with $p^{Y^h}_1(x)=p^Y_1(x)+h(p^A_1(x)-p^B_1(x))$.
In (\ref{185}), the estimate involving the density $f(N,Y)$ is similar, while the estimate for
$\E^{Y^h}_0[\check Z^{z,{\rm pin}}_{N,Y^h}]$ reduces to estimating
$$
\begin{aligned}
G^{X-Y}-G^{X-Y^h} &= \frac{1}{(2\pi)^d}\int_{[-\pi,\pi]^d}
\Big(\frac{1}{1-\phi_X(k)\overline\phi_Y(k)}-\frac{1}{1-\phi_X(k)\overline\phi_{Y^h}(k)}\Big)dk \\
&= \frac{h}{(2\pi)^d}\int_{[-\pi,\pi]^d}
\frac{\phi_X(\overline\phi_B-\overline\phi_A)}{\big(1-\phi_X\overline\phi_Y\big)\big(1-\phi_X\overline\phi_{Y^h}\big)}dk,
\end{aligned}
$$
where $\phi_X(k)=\sum_{x}e^{ik\cdot x}p^X_1(x)$, $\overline\phi_X(k)=\phi_X(-k)$, and $\overline\phi_Y(k)$, $\overline\phi_A(k)$
and $\overline\phi_B(k)$ are defined similarly. Note that in $d\geq 4$, $\int\big|\frac{\phi_X(\overline\phi_B-\overline\phi_A)}{(1-\phi_X\overline\phi_Y)^2}\big|dk<\infty$.
Therefore based on Taylor expansion in $h$, all calculations carry through as long as
\be\label{Qcond}
Q:=\int\frac{\phi_X(\overline\phi_B-\overline\phi_A)}{(1-\phi_X\overline\phi_Y)^2}dk\neq0
\ee
and $h$ is chosen to have the same sign. When $X$ and $Y$ are simple random walks, we have $Q=0$ for any choice of
$A$, $B$, $p^A_1$ and $p^B_1$ due to symmetry. In particular, changing the drift for the simple random walk fails.
On the other hand, if $S$ contains enough points so as to break symmetry, then it is reasonable to expect the existence of
$A$, $B$, $p^A_1$ and $p^B_1$ which give $Q\neq 0$.
When such $A, B, p^A_1$ and $p^B_1$ exist, we may even take $A$ and $B$ to be singletons in $S$. We were not able to verify
(\ref{Qcond}) for some $A,B\subset S$ for general random walks, such as for all walks with zero mean and finite variance and whose
support $S$ contains at least two points which are not related by reflections or permutations of coordinates. However,
when $X$ and $Y$ are i.i.d.\ so that $\phi_X=\phi_Y$, $\phi_X\geq 0$, and $0\in S$, it is easily seen that $Q>0$ for $B=\{0\}$ and
$p^A_1 = p^X_1$. This includes random walks $X$ which are symmetric with $p^X_1(0)\geq \frac{1}{2}$, as well as walks $X$ that
can be expressed as the difference of two i.i.d.\ random walks.

\subsection{Proof of Theorem \ref{T:cptst} in discrete time: $d=4$}
For $d=4$, in the representation (\ref{diszfree}), we have $K(n)=p^{X-Y}_n(0)/G^{X-Y} \sim C n^{-2}$ which has
infinite first moment. Thus $d=4$ corresponds to the case $\alpha=1$ in \cite{DGLT07} for the pinning model. In
\cite{DGLT07}, the case $\alpha=1$ was left out. However, as we will show below, there is no difficulty in extending
the fractional moment method to the $d=4$ case, and we expect the same to be true for the $\alpha=1$ case for the
pinning model.

As in $d\geq 5$, it suffices to verify (\ref{rhodis}). What differs in $d=4$ is that
$\sum_{i=R}^\infty i^{1-\frac{d\gamma}{2}}= \sum_{i=R}^\infty i^{1-2\gamma}=\infty$ for any $\gamma\in (0,1)$
and $R\in\N$. Hence a more careful estimate of $\check \varrho$ than in (\ref{456}) is needed. By Theorem 2.1 of
\cite{G07} and super-multiplicativity of $(\E^Y_0[\check Z^{z, \rm pin}_{n,Y}])_{n\in\N}$, we have
$\E^Y_0[\check Z^{z, \rm pin}_{N,Y}] \leq e^{C N (z-1)}$ for some $C>0$ uniformly in $z>1$ sufficiently close to 1
and $N\in\N$. Therefore the same choice $L=(z-1)^{-1}$ as in $d\geq 5$ ensures that
$\sup_{1\leq i\leq L}\check A^{z, \rm pin}_i \leq C<\infty$ uniformly for $z>1$ close to 1. Fix $\eps>0$ small, then
let $\gamma \in (0,1)$ such that $2\gamma-1>1-\eps$. Analogous to (\ref{456}), we have
\be\label{rhodec4}
\check \varrho \leq \sum_{i=0}^{L^{1-\eps}} \frac{C}{(L-i)^{2\gamma-1}}
+ \sum_{i=L^{1-\eps}}^{L-1} \frac{C\check A^{z, \rm pin}_i}{(L-i)^{2\gamma-1}}
\leq \frac{CL^{1-\eps}}{(L-L^{1-\eps})^{2\gamma-1}} + C L^{2-2\gamma} \max_{L^{1-\eps}\leq i\leq L}\check A^{z, \rm pin}_i.
\ee
Therefore to show $\check \varrho<1$ for some $z>1$, it suffices to show that with $L=(z-1)^{-1}$,
\be\label{424}
\lim_{z\downarrow 1} L^{2-2\gamma}\!\!\!\!  \max_{L^{1-\eps}\leq N \leq L} \check A^{z, \rm pin}_N = 0.
\ee
Tracing through the arguments for $d\geq 5$, we see that analogous to (\ref{462}), for $h=1/\sqrt{L}=\sqrt{z-1}$,
uniformly for $L^{1-\eps}\leq N\leq L$ and $z>1$ sufficiently close to 1, we have
\be
\check A^{z, \rm pin}_N \leq C \E^{Y^h}_0[\check Z^{z, \rm pin}_{N, Y^h}]^\gamma
\leq C \E^{K_*}\left[\exp\left\{-c\sqrt{z-1}\, \big|\iota^*\cap [1,(z-1)^{\eps-1}]\big|\right\}\right]^\gamma,
\ee
where $\iota^*$ is a renewal process on $\N_0$ with inter-arrival probability distribution $K_*$ satisfying the
property that $K_*(n) \sim C n^{-2}$ for some $C>0$. Set $M=(z-1)^{\eps-1}$. Then
\be\nonumber
0\leq \lim_{z\downarrow 1} L^{2-2\gamma}\!\!\!\!  \max_{L^{1-\eps}\leq N \leq L} \check A^{z, \rm pin}_N
\leq \lim_{M\to\infty} C M^{\frac{2-2\gamma}{1-\eps}} \E^{K_*}\left[\exp\left\{-c M^{-\frac{1}{2(1-\eps)}}\,\big|\iota^*\cap [1,M]\big|\right\}\right]^\gamma
=0,
\ee
where we applied Proposition \ref{P:rewest} with $\delta_1 = \frac{1}{2(1-\eps)}$ and $1-\delta_2 =
\frac{2-2\gamma}{\gamma(1-\eps)}$, which satisfy the condition $0<\delta_1<\delta_2<1$ if $\eps>0$ is small, and
$\gamma\in(0,1)$ is then chosen sufficiently close to 1.
\qed

\section{Gap between critical points: continuous time}

\subsection{Proof of Theorem \ref{T:cptst} in continuous time: $d\geq 5$}
As in discrete time, we split the proof into three parts: representation for $Z^\beta_{t,Y}$ and
$Z^{\beta,\rm pin}_{t,Y}$; fractional moment method; change of measure. Compared to the discrete time case,
the main complication here is to suitably discretize time so that the fractional moment inequality
(\ref{fracineq}) can be applied. The change of measure argument however becomes much simpler.

\bigskip
\noindent {\bf Representation for $Z^\beta_{t,Y}$ and  $Z^{\beta, \rm pin}_{t,Y}$.} We now Taylor expand
$e^{\beta L_t(X,Y)}$. Let $p_s(\cdot)$ be the transition probability kernel of a rate 1 continuous time simple
random walk on $\Z^d$. Let $G_{1+\rho}=\int_0^\infty p_{(1+\rho)s}(0)ds$, $K_{1+\rho}(s) = p_{(1+\rho)s}(0)/G_{1+\rho}$,
$\bar \beta = \beta G_{1+\rho}$, and $\bar Z^{\bar\beta}_{t,Y} = Z^{\beta}_{t,Y}$. Then

\begin{eqnarray}
\bar Z^{\bar \beta}_{t,Y}
&=&
\E^X_0\left[1+\sum_{m=1}^\infty \frac{\beta^m}{m!}\Big(\int_0^t 1_{\{X_s=Y_s\}}ds\Big)^m\right]  \nonumber\\
&=&
\E^X_0\left[1+\sum_{m=1}^\infty \beta^m \idotsint\limits_{0<s_1\cdots<s_m<t}
1_{\{X_{s_1}=Y_{s_1},\cdots, X_{s_m}=Y_{s_m}\}}ds_1\cdots ds_m\right]  \nonumber \\
&=&
1 + \sum_{m=1}^\infty \beta^m \idotsint\limits_{0<s_1\cdots<s_m<t}
p_{s_1}(Y_{s_1})p_{s_2-s_1}(Y_{s_2}-Y_{s_1})\cdots p_{s_m-s_{m-1}}(Y_{s_m}-Y_{s_{m-1}})ds_1\cdots ds_m \nonumber \\
&=&
1 + \sum_{m=1}^\infty
\idotsint\limits_{ \atop s_0=0<s_1\cdots<s_m<t}
\prod\limits_{i=1}^{m} \Big(K_{1+\rho}(s_i-s_{i-1})w(\bar\beta, s_i-s_{i-1}, Y_{s_i}-Y_{s_{i-1}})\Big) ds_1\cdots
ds_m, \label{rwpartrepf}
\end{eqnarray}
where
\be\label{wcont}
w(\bar\beta, s_i-s_{i-1}, Y_{s_i}-Y_{s_{i-1}}) = \frac{\bar\beta\, p_{s_i-s_{i-1}}(Y_{s_i}-Y_{s_{i-1}})}
{p_{(1+\rho)(s_i-s_{i-1})}(0)}.
\ee
If we denote $\bar Z^{\bar\beta, \rm pin}_{t,Y} = \beta Z^{\beta, \rm pin}_{t,Y}$, then similarly,
\be\label{rwpartrep}
\bar Z^{\bar \beta,\rm pin}_{t,Y}= K_{1+\rho}(t)w(\bar\beta,t, Y_t)\! + \!\!\sum_{m=1}^\infty
\!\!\!\!\!\!\!\!\!\!\!\!\!
\idotsint\limits_{ \atop s_0=0<s_1\cdots<s_m<s_{m+1}=t}
\!\!\!\!\!\!\!\!\!\!\!\!\!
\prod\limits_{i=1}^{m+1} K_{1+\rho}(s_i\!-\!s_{i-1})w(\bar\beta, s_i\!-\!s_{i-1}, Y_{s_i}\!-\!Y_{s_{i-1}}) ds_1\cdots
ds_m.
\ee
Note that (\ref{rwpartrep}) casts $\bar Z^{\bar\beta,\rm pin}_{t,Y}$ in the same form as (\ref{disz}),
except that the underlying renewal process is in continuous time with return time distribution $K_{1+\rho}(s)ds$.
Since
\be\label{sth}
\E^Y_0[w(\bar\beta, s_i-s_{i-1}, Y_{s_i}-Y_{s_{i-1}})]=\bar\beta,
\ee
and $K_{1+\rho}(\cdot)$ defines a recurrent
renewal process on $[0,\infty)$, $\E^Y_0[\bar Z^{\bar\beta, \rm pin}_{t,Y}]$ is the partition function of
a homogeneous pinning model (in continuous time) with critical point $\bar\beta^{\rm ann}_{\rm c}=1$, or
equivalently,
\be\label{abccts}
\beta^{\rm ann}_{\rm c} = \frac{\bar\beta^{\rm ann}_{\rm c}}{G_{1+\rho}} = \frac{1}{G_{1+\rho}}.
\ee

\noindent
{\bf Fractional moment method.}
Analogous to (\ref{checkZsplit}), for fixed $L\in\N$, we have the decomposition
\be
\bar Z^{\bar \beta}_{t,Y} =
\bar Z^{\bar\beta}_{L,Y} + \iint\limits_{0\leq u<L<v\leq t}\!\!\!\!\!\!\! K_{1+\rho}(v-u)w(\bar\beta, v-u,Y_v-Y_u)
\bar Z^{\bar\beta, \rm pin}_{u,Y} \bar Z^{\bar\beta}_{t-v, \theta_vY} (1+\delta_0(u))dudv,
\label{rwpartrep2}
\ee
where $\theta_vY = (Y_{v+s}-Y_v)_{s\geq 0}$ denotes a shift in $Y$, $\delta_0(u)$ is the delta function at $0$, and
$\bar Z^{\bar\beta, \rm pin}_{0,Y}=1$. In the continuous setting, the analogue of (\ref{fracineq}),
$(\int |a(x)|dx)^\gamma \leq \int |a(x)|^\gamma dx$ for $\gamma\in(0,1)$, is false in general. Therefore we
need to discretize the integrals in (\ref{rwpartrep2}). In order to obtain uniform control for the integrand in
(\ref{rwpartrep2}) on intervals, it turns out to be more suitable to study the following quantities in place
of $\bar Z^{\bar\beta}_{t,Y}$ and $\bar Z^{\bar\beta, \rm pin}_{t,Y}$.
\be\label{Ztylr}
\begin{aligned}
\bar Z^{\bar\beta, 1}_{t,Y}\!\! &=& 1\quad+\quad \sum_{m=1}^\infty\!\!\!\!\!
 \idotsint\limits_{\atop s_0=0<s_1\cdots<s_m<t}
\!\!\!\!\!
\prod_{i=1}^{m}K_{1+\rho}(s_i-s_{i-1}) \prod\limits_{i=2}^m w(\bar\beta,s_i-s_{i-1}, Y_{s_i}-Y_{s_{i-1}}) ds_1\cdots
ds_m, \\
\bar Z^{\bar \beta, {\rm pin1}}_{t,Y} \!\! &=&\!\!\! K_{1+\rho}(t) + \!
\sum_{m=1}^\infty \!\!\!\!\!\!\!\!\!\!\!\!\!
\idotsint\limits_{\atop s_0=0<s_1\cdots<s_m<s_{m+1}=t} \!\!\!\!\!\!\!\!\!\!\!\!
\prod_{i=1}^{m+1}K_{1+\rho}(s_i\!-\!s_{i-1})\!\!
 \prod\limits_{i=2}^{m+1} w(\bar\beta,s_i\!-\!s_{i-1}, Y_{s_i}\!-\!Y_{s_{i-1}}) ds_1\cdots ds_m,\\
\bar Z^{\bar \beta, {\rm pin2}}_{t,Y}
\!\! &=&\!\!\!\!\!\!\!\!\!
K_{1+\rho}(t) + \!\!
\sum_{m=1}^\infty \!\!\!\!\!\!\!\!\!\!\!\!\!
\idotsint\limits_{\atop s_0=0<s_1\cdots<s_m<s_{m+1}=t} \!\!\!\!\!\!\!\!\!\!\!\!
\prod_{i=1}^{m+1}K_{1+\rho}(s_i\!-\!s_{i-1}) \!
\prod\limits_{i=2}^{m} w(\bar\beta,s_i\!-\!s_{i-1}, Y_{s_i}\!-\!Y_{s_{i-1}}) ds_1\cdots ds_m,\\
\end{aligned}
\ee
where $\prod_{i=2}^m w=1$ if $m=1$. Note that $\bar Z^{\bar\beta, 1}_{t,Y}$ differs from
$\bar Z^{\bar\beta}_{t,Y}$ in that the factor $w(\bar\beta, s_1, Y_{s_1})$ in (\ref{rwpartrepf})
has been omitted, while $\bar Z^{\bar\beta, \rm pin1}_{t,Y}$ (resp.\ $\bar Z^{\bar\beta, \rm pin2}_{t,Y}$)
differs from $\bar Z^{\bar\beta, \rm pin}_{t,Y}$ in that the factors $w(\bar\beta, t, Y_t)$ and
$w(\bar\beta, s_1, Y_{s_1})$ (resp.\ as well as $w(\bar\beta, t-s_m, Y_t-Y_{s_m})$) in (\ref{rwpartrep})
have been omitted. Omitting these random factors will provide flexibility in adjusting the lengths of
the renewal gaps $(s_i-s_{i-1})_{i\in\N}$.

Note that
\be\label{Wbound}
w(\bar\beta, v-u, Y_v-Y_u)=\frac{\bar\beta\, p_{v-u}(Y_v-Y_u)}{p_{(1+\rho)(v-u)}(0)}\leq
\frac{\bar\beta\, p_{v-u}(0)}{p_{(1+\rho)(v-u)}(0)}\leq C
\ee
for some $C\in(1,\infty)$ independent of $v-u\geq 0$ and $\bar\beta \in [1,2]$, which is furthermore uniformly
bounded for $\rho\in [0,1]$. Therefore
\be\label{Ztycomp}
\bar Z^{\bar\beta}_{t,Y} \leq C \bar Z^{\bar\beta, 1}_{t,Y}.
\ee
By the monotonicity of $Z^\beta_{t,Y}=\bar Z^{\bar\beta}_{t,Y}$ in $t$, to show
$\beta<\beta^*_{\rm c}$ (i.e., $\sup_{t\geq 0} Z^{\beta}_{t,Y}<\infty$ a.s.\ w.r.t.\ $Y$), it suffices to show
that for $\bar\beta = \beta G_{1+\rho}$, there exists $\gamma\in (0,1)$
such that
\be\label{579}
\sup_{t\geq 0} \E^Y_0\big[\big(\bar Z^{\bar\beta, 1}_{t,Y}\big)^\gamma\big] < \infty.
\ee
Note that $\bar Z^{\bar\beta, 1}_{t,Y}$ is increasing in $t$ for every $Y$, therefore we may assume
$t\in\N$. Similar to (\ref{rwpartrep2}), we
have
\begin{eqnarray}
\bar Z^{\bar\beta, 1}_{t,Y} &=& \bar Z^{\bar\beta, 1}_{L,Y} +
\int_L^t K_{1+\rho}(v)\bar Z^{\bar\beta}_{t-v, \theta_vY}dv + \!\!\!\!\!\!
\iint\limits_{0< u<L<v< t}\!\!\!\!\!\!\! K_{1+\rho}(v-u)w(\bar\beta, v-u,Y_v-Y_u)
\bar Z^{\bar\beta, \rm pin1}_{u,Y} \bar Z^{\bar\beta}_{t-v, \theta_vY} dudv \nonumber \\
&=& \bar Z^{\bar\beta, 1}_{L,Y} +
\sum_{j=L}^{t-1}\int\limits_{j}^{j+1}K_{1+\rho}(v) \bar Z^{\bar\beta}_{t-v, \theta_v Y}dv \label{Ztydecomp} \\
&& \qquad \qquad \qquad + \sum_{i=0}^{L-1} \sum_{j=L}^{t-1}\ \nonumber
\iint\limits_{i<u<i+1\atop j<v<j+1}
K_{1+\rho}(v-u) w(\bar\beta, v-u, Y_v-Y_u) \bar Z^{\bar\beta, \rm pin1}_{u,Y} \bar Z^{\bar\beta}_{t-v,\theta_vY}dudv.
\end{eqnarray}
We will establish uniform estimates on the integrand for each integral in (\ref{Ztydecomp}) by bounding
$\bar Z^{\bar\beta}_{t-v, \theta_vY}$ in terms of $\bar Z^{\bar\beta, 1}_{t-j-1, \theta_{j+1}Y}$ and
bounding $\bar Z^{\bar\beta,\rm pin1}_{u,Y}$ in terms of $\bar Z^{\bar\beta, \rm pin2}_{i,Y}$.

We first make a few observations which will come in handy. Note that for all $s\in [0,1]$ and all realizations of $Y$,
\be
\begin{aligned}
\bar Z^{\bar\beta}_{s,Y}  &= Z^{\beta}_{s,Y} = \E^X_0[e^{\beta L_s(X,Y)}]\leq e^\beta, \\
\bar Z^{\bar\beta, \rm pin}_{s,Y} &= \beta Z^{\beta, \rm pin}_{s,Y} = \beta\, \E^X_0[e^{\beta L_s(X,Y)}1_{\{X_s=Y_s\}}]
\leq \beta e^\beta.
\end{aligned}
\ee
Next note that
\be\label{Kbound}
C_{\rho} = \sup_{u\geq 0 \atop 0\leq s\leq 1} \frac{K_{1+\rho}(u)}{K_{1+\rho}(u+s)}<\infty,
\ee
which is uniformly bounded for $\rho\in [0,1]$.

If $v\in (j, j+1)$ for some $L\leq j\leq t-1$, then by the same decomposition as (\ref{rwpartrep2})
with $s_1,s_2, j+1$ now playing the roles of $u,v,L$ and by the observations above, we have
\begin{eqnarray}
\bar Z^{\bar\beta}_{t-v, \theta_vY}
\!\!\!\!\! &=&\!\!\!\!\!
\bar Z^{\bar\beta}_{j+1-v, \theta_vY} +     \nonumber
\!\!\! \iint\limits_{v\leq s_1<j+1 \atop j+1<s_2<t} \!\!\!
K_{1+\rho}(s_2-s_1)w(\bar\beta, s_2\!-\!s_1, Y_{s_2}\!-\!Y_{s_1}) \bar
Z^{\bar\beta, \rm pin}_{s_1-v, \theta_vY} \bar Z^{\bar\beta}_{t-s_2, \theta_{s_2}Y} (1\!+\!\delta_v(s_1))ds_1ds_2 \\
&\leq&\!\!\!\!\!
C + C \int_{j+1}^t K_{1+\rho}(s_2-j-1)  \bar Z^{\bar\beta}_{t-s_2, \theta_{s_2}Y} ds_2
= C\bar Z^{\bar\beta,1}_{t-j-1,\theta_{j+1}Y}, \label{582}
\end{eqnarray}
where $C<\infty$ is independent of $t, v, Y$, $\bar\beta\in [1,2]$, and furthermore is uniformly bounded for
$\rho\in[0,1]$.

If $u\in (i,i+1)$ for some $0\leq i\leq L-1$, then by a similar decomposition as above, we have
\begin{eqnarray}\nonumber
\bar Z^{\bar\beta, {\rm pin1}}_{u,Y} &=&
\int\limits_{i<s_2\leq u} K_{1+\rho}(s_2) \bar Z^{\bar\beta, \rm pin}_{u-s_2,\theta_{s_2}Y} (1+\delta_u(s_2))ds_2 \\
&& \quad +
\iint\limits_{0<s_1<i<s_2\leq u} K_{1+\rho}(s_2-s_1) w(\bar\beta, s_2-s_1,Y_{s_2}-Y_{s_1})
\bar Z^{\bar\beta, \rm pin1}_{s_1,Y}
\bar Z^{\bar\beta, \rm pin}_{u-s_2,\theta_{s_2}Y}(1+\delta_u(s_2))ds_1ds_2 \nonumber \\
&\leq& C K_{1+\rho}(i) + C\int\limits_{0<s_1<i} K_{1+\rho}(i-s_1)\bar Z^{\bar\beta, \rm pin1}_{s_1,Y}ds_1
\ =\ C \bar Z^{\bar\beta, \rm pin2}_{i,Y}. \label{583}
\end{eqnarray}
Substituting the bounds (\ref{Wbound}), (\ref{Kbound})--(\ref{583}) into (\ref{Ztydecomp}) gives
\begin{eqnarray}
\bar Z^{\bar\beta, 1}_{t,Y}\!\!\!\! &\leq&\!\!\!\!
\bar Z^{\bar\beta, 1}_{L,Y} + C' \sum_{j=L}^{t-1} K_{1+\rho}(j+1) \bar Z^{\bar\beta, 1}_{t-j-1, \theta_{j+1}Y}
+ C' \sum_{i=0}^{L-1} \sum_{j=L}^{t-1} K_{1+\rho}(j+1-i) \bar Z^{\bar\beta, \rm pin2}_{i,Y}
\bar Z^{\bar\beta, 1}_{t-j-1,\theta_{j+1}Y} \nonumber \\
&\leq&
\bar Z^{\bar\beta, 1}_{L,Y} +
C \sum_{i=0}^{L-1} \sum_{j=L}^{t-1} K_{1+\rho}(j+1-i) \bar Z^{\bar\beta, \rm pin2}_{i,Y}
\bar Z^{\bar\beta, 1}_{t-j-1,\theta_{j+1}Y}, \label{tildeZbd3}
\end{eqnarray}
where $C<\infty$ is independent of $t$, $Y$, $\bar\beta\in[1,2]$, and can be chosen uniformly for
$\rho\in[0,1]$.

Fix $\gamma\in (0,1)$ such that $\frac{d\gamma}{2}>2$ for $d\geq 5$. Denote
$\bar A^{\bar\beta, 1}_t=\E^Y_0\big[\big(\bar Z^{\bar\beta, 1}_{t,Y}\big)^\gamma\big]$ and
$\bar A^{\bar\beta, \rm pin2}_t= \E^Y_0\big[\big(\bar Z^{\bar\beta, {\rm pin2}}_{t,Y}\big)^\gamma\big]$.
Then the same calculations as those leading to (\ref{453}) yields
\be\label{rhocts}
\bar A^{\bar\beta,1}_t \leq \bar A^{\bar\beta, 1}_{L} +
\varrho \sup_{0\leq j\leq t-L}\bar A^{\bar\beta, 1}_j \qquad \mbox{with}\qquad
\varrho = C\left(\sum_{i=0}^{L-1} \frac{\bar A^{\bar\beta,\rm pin2}_i}{(L-i)^{\frac{d\gamma}{2}-1}}\right),
\ee
where $C<\infty$ is independent of $t$ and $\bar\beta\in [1,2]$, and can be chosen uniformly for
$\rho\in [0,1]$. As in the discrete time case, we aim to show $\varrho<1$.

Note that $\bar A^{\bar\beta, \rm pin2}_s \leq \E^Y_0[\bar Z^{\bar\beta, \rm pin2}_{s,Y}]^\gamma\leq \E^Y_0[\bar
Z^{\bar\beta,\rm pin}_{s,Y}]^\gamma \leq \E^Y_0[\bar Z^{\bar\beta}_{s,Y}]^\gamma$ by Jensen and (\ref{sth}),
where we see from (\ref{rwpartrepf}) that
$\E^Y_0[\bar Z^{\bar\beta}_{s,Y}]$ is the partition function of a continuous time homogeneous pinning model
with return time distribution $K_{1+\rho}(\cdot)$ and critical point $\bar\beta^{\rm ann}_{\rm c}=1$. For $d\geq 5$,
it is easy to verify (by law of large numbers and elementary large deviation estimates for the number of returns
of the renewal process before time $s$) that
\be\label{517}
\E^Y_0[\bar Z^{\bar\beta}_{s,Y}]\leq Ce^{C(\bar\beta-1)s}
\ee
for some $C\in (0,\infty)$ independent of $s\geq 0$ and $\bar\beta\in [1,2]$, and is furthermore uniformly bounded for
$\rho\in [0,1]$. As in the discrete time case, we choose
\be\label{Lcts}
L=(\bar\beta-1)^{-1}.
\ee
In view of (\ref{579}) and (\ref{rhocts}), and by the same arguments as those leading to (\ref{ALbound}) in the discrete
time case, to show $\beta^*_{\rm c}>\beta^{\rm ann}_{\rm c}$ for any $\rho>0$, it suffices to show that
\be\label{supAL}
\lim_{\bar\beta \downarrow 1} \sup_{L-R\leq t\leq L}\bar A^{\bar\beta, \rm pin2}_t =0,
\ee
where $R\in\N$ is large and fixed and can be chosen uniformly for $\rho\in [0,1]$.
On the other hand, showing
\be\label{sharpcptbd}
\beta^*_{\rm c}-\beta^{\rm ann}_{\rm c} \geq a\rho
\ee
for some $a>0$ and all $\rho\in [0,1]$ reduces to showing that: (1) the convergence in (\ref{supAL}) is in fact uniform for
$\rho \in [\rho_0, 1]$ for any $0<\rho_0\leq 1$, which implies that
$\inf_{\rho\in [\rho_0, 1]}(\bar\beta_{\rm c}^*-1) >0$ where $\bar\beta^*_{\rm c}=G_{1+\rho}\beta^*_{\rm c}$,
and hence $\inf_{\rho\in [\rho_0,1]} (\beta^*_{\rm c}-\beta^{\rm ann}_{\rm c})>0$;
(2) for $\bar\beta = 1+a\rho$ with $a>0$ sufficiently small, $L=(\bar\beta-1)^{-1}$, and $R\in\N$ large and independent
of $\rho\in [0,1]$,
\be\label{supAL2}
\limsup_{\rho\downarrow 0} \sup_{L-R\leq t\leq L} \bar A^{\bar\beta, \rm pin2}_t < 1,
\ee
which implies that for some $\rho_0\in (0,1]$, $\bar\beta_{\rm c}^*-1 = G_{1+\rho}(\beta^*_{\rm c}-\beta^{\rm ann}_{\rm c})\geq a\rho$
for all $\rho\in [0, \rho_0]$.
\bigskip

\noindent
{\bf Change of measure.} We now prove (\ref{supAL}) and (\ref{supAL2}), where the convergence in (\ref{supAL}) will
be shown to be uniform in $\rho\in [\rho_0, 1]$ for any $0<\rho_0\leq 1$. Here, the appropriate
change of measure for the disorder $Y$ is simply to increase the jump rate of the random walk $Y$. Let $Y^{\rho+h}$
be a simple random walk on $\Z^d$ with jump rate $\rho+h$ for some $h>0$, then the path measures $(Y_s)_{0\leq s\leq t}$
and $(Y^{\rho+h}_s)_{0\leq s\leq t}$ are equivalent, and the Radon-Nikodym derivative of the
law of $(Y^{\rho+h}_s)_{0\leq s\leq t}$ w.r.t.\ that of $(Y_s)_{0\leq s\leq t}$ is given by
\be\nonumber
f(t,Y)= e^{-h t}(1+h\rho^{-1})^{N_t(Y)},
\ee
where $N_t(Y)$ is the number of jumps of $Y$ in $[0,t]$. Then as in (\ref{185}),
\be
\bar A^{\bar\beta, \rm pin2}_t
= \E^{Y^{\rho+h}}_0\big[ f(t,Y^{\rho+h})^{-1}
\big(\bar Z^{\bar\beta, {\rm pin2}}_{t,Y^{\rho+h}}\big)^\gamma\big] \leq
\E^{Y}_0\big[f(t,Y)^{-\frac{\gamma}{1-\gamma}}]^{1-\gamma}\
\E^{Y^{\rho+h}}_0\big[\bar Z^{\bar\beta, {\rm pin2}}_{t,Y^{\rho+h}}\big]^\gamma. \label{163}
\ee
Note that
\begin{eqnarray}
\E^{Y}_0\big[f(t,Y)^{-\frac{\gamma}{1-\gamma}}] &=& e^{\frac{\gamma h t}{1-\gamma}}
\E^Y_0\big[(1+h\rho^{-1})^{-\frac{\gamma N_t}{1-\gamma}}\big] =
e^{\frac{\gamma h t}{1-\gamma}} \sum_{n=0}^\infty e^{-\rho t} \frac{(\rho t)^n}{n!}
(1+h\rho^{-1})^{-\frac{\gamma n}{1-\gamma}} \nonumber \\
&=& \exp\left\{\Big(\rho(1+h \rho^{-1})^{-\frac{\gamma}{1-\gamma}}-\rho +
  \frac{\gamma h}{1-\gamma}\Big)t\right\}
\leq \exp\left\{\frac{\gamma h^2 t}{2\rho(1-\gamma)^2}\right\}, \label{densitypower}
\end{eqnarray}
where second order Taylor expansion in $h$ in the exponent provides a true upper bound.
For $L-R\leq t\leq L$, if we choose $h=\frac{\sqrt\rho}{\sqrt L}$, then the first term in (\ref{163})
is bounded and independent of $\rho$, $\bar\beta$ and $t$. Thus it only remains to estimate
$\E^{Y^{\rho+h}}_0\big[\bar Z^{\bar\beta, {\rm pin2}}_{t,Y^{\rho+h}}\big]$.

First note that $\E^{Y^{\rho+h}}_0\big[\bar Z^{\bar\beta, {\rm pin2}}_{t,Y^{\rho+h}}\big]\leq C
\E^{Y^{\rho+h}}_0\big[\bar Z^{\bar\beta, \rm pin}_{t,Y^{\rho+h}}\big]$ for some $C>0$ independent of $\rho\geq 0$,
$\bar\beta\in [1,2]$ and $t\geq 0$, because each term in the expansion for $\bar Z^{\bar\beta, \rm pin}_{t,Y}$ in
(\ref{rwpartrep}) differs from the corresponding term in (\ref{Ztylr}) for $\bar Z^{\bar\beta, {\rm pin2}}_{t,Y}$
by at most two factors of $w$, and
$\E^{Y^{\rho+h}}_0[w(\bar\beta, v-u,Y_v-Y_u)] = \frac{\bar\beta\, p_{(1+\rho+h)(v-u)}(0)}{p_{(1+\rho)(v-u)}(0)}\geq C$
for some $C>0$ independent of $\rho\geq 0$, $h\in[0,1]$, $\bar\beta\in [1,2]$ and $v-u\geq 0$.
Recall $G_{1+\rho}=\int_0^\infty p_{(1+\rho)s}(0)ds$,
\begin{eqnarray}
&& \E^{Y^{\rho+h}}_0\Big[\bar Z^{\bar\beta, \rm pin}_{t,Y^{\rho+h}}\Big] \nonumber\\
&=&  \left(\frac{\bar\beta}{G_{1+\rho}}\right) p_{(1+\rho+h)t}(0) + \sum_{m=1}^\infty\
\idotsint\limits_{0=s_0<s_1\cdots<s_m<s_{m+1}=t} \left(\frac{\bar\beta}{G_{1+\rho}}\right)^{m+1}
\prod_{i=1}^{m+1} p_{(1+\rho+h)(s_i-s_{i-1})}(0)\, ds_1\cdots ds_m \nonumber \\
&=& \frac{(1+\rho)\bar\beta}{1+\rho+h} K_{1+\rho+h}(t) + \!\!\!\!\!\!\!\!\!
\idotsint\limits_{\atop 0=s_0<s_1\cdots<s_m<s_{m+1}=t} \!\!\!\!\!\! \left(\frac{(1+\rho)\bar\beta}{1+\rho+h}\right)^{m+1}
\prod_{i=1}^{m+1} K_{1+\rho+h}(s_i-s_{i-1})\, ds_1\cdots ds_m, \label{Kmod}
\end{eqnarray}
where $K_{1+\rho+h}(s) = p_{(1+\rho+h)s}(0)/G_{1+\rho+h}$ with $G_{1+\rho+h} = \int_0^\infty p_{(1+\rho+h)s}(0)ds =
\frac{(1+\rho)G_{1+\rho}}{1+\rho+h}$.

Denote $\bar\beta' = \frac{(1+\rho)\bar\beta}{1+\rho+h}$. Let $\sigma^{\rho+h}=(0, \sigma^{\rho+h}_1,\sigma^{\rho+h}_2,
\cdots)$ be a renewal sequence on $[0,\infty)$ with inter-arrival law $K_{1+\rho+h}(\cdot)$, and let
$\E^{K_{1+\rho+h}}[\cdot]$ denote expectation w.r.t.\ $\sigma^{\rho+h}$. Then in view of (\ref{Kmod}),
\be
\E^{K_{1+\rho+h}}\left[(\bar\beta')^{1+|\sigma^{\rho+h} \cap [0,t]|}\ 1_{\{\sigma^{\rho+h} \cap [t,t+1]\neq\emptyset\}}\right]
\geq \inf_{u\geq 0,\atop 0\leq s\leq 1} \frac{K_{1+\rho+h}(u+s)}{K_{1+\rho+h}(u)}
\ \E^{Y^{\rho+h}}_0\big[\bar Z^{\bar\beta, \rm pin}_{t,Y^{\rho+h}}\big]. \nonumber
\ee
Recall the definition of $C_{1+\rho}$ from (\ref{Kbound}), we then have
\be\label{168}
\E^{Y^{\rho+h}}_0\big[\bar Z^{\bar\beta, \rm pin}_{t,Y^{\rho+h}}\big] \leq C_{\rho+h}
\E^{K_{1+\rho+h}}\left[(\bar\beta')^{1+|\sigma^{\rho+h} \cap [0,t]|}\right].
\ee

Now to prove (\ref{supAL2}), we recall that $L=(\bar\beta-1)^{-1}$ and hence
$h=\frac{\sqrt \rho}{\sqrt L}=\sqrt{\rho(\bar\beta-1)}$.
Therefore there exists $\bar\beta_0>1$ sufficiently small such that for all $\rho>0$ and $\bar\beta \in [1, \bar\beta_0]$,
\be
\bar\beta' = \frac{(1+\rho)\bar\beta}{1+\rho+h} \leq (1+\bar\beta -1)\Big(1-\frac{\sqrt{\rho(\bar\beta-1)}}{2(1+\rho)}\Big).
\ee
First note that by our choice $\bar\beta = 1+a\rho$, we have $\bar\beta' \leq 1-\rho\sqrt{a}/8$ for all $\rho\in [0,1]$
if $0<a<1/64$. Next note that $C_{\rho+h}$ is uniformly bounded for $\rho\in [0,1]$ and $\bar\beta\in [1,2]$. For $d\geq 5$,
by the local central limit theorem, there exists an inter-arrival probability distribution $K_*$ on $(0,\infty)$ with finite
first moment $m=\int_0^\infty s K_*(s)ds$, such that $K_*$  stochastically dominates $K_{1+\rho+h}$ for all $h\in [0,1]$ and
$\rho\in [0, 1]$. Namely, $\int_t^\infty K_*(s)ds \geq \int_t^\infty K_{1+\rho+h}(s)ds$ for all $t\geq 0$, $h\in[0,1]$ and
$\rho\in [0, 1]$. Combining the above observations, we have
\be\label{596}
\limsup_{\rho\downarrow 0}\!\!\! \sup_{L-R\leq t\leq L} \!\!\!\! \bar A^{\bar\beta, \rm pin2}_t
\! \leq\!
C \limsup_{\rho\downarrow 0}\!\!\! \sup_{L-R\leq t\leq L} \!\!\!\! \E^{Y^{\rho+h}}_0\big[\bar Z^{\bar\beta, \rm pin}_{t,Y^{\rho+h}}\big]
\!\leq\!
C \limsup_{\rho\downarrow 0}\E^{K_*}[(1-\rho\sqrt{a}/8)^{|\iota^*\cap [0,L-R]|}],
\ee
where $\iota^*$ is a renewal process on $[0,\infty)$ with return time distribution $K_*$. By the law of large numbers, a.s.\ w.r.t.\
$\iota^*$,
\be\nonumber
\lim_{\rho\downarrow 0} (1-\rho\sqrt{a}/8)^{|\iota^*\cap [0,L-R]|}
= \lim_{\rho\downarrow 0} \exp\left\{-\frac{\rho \sqrt{a}}{8}\cdot \frac{(a\rho)^{-1}-R}{m}\right\} = \exp\left\{-\frac{1}{8m\sqrt{a}}\right\},
\ee
which can be made arbitrarily small if $a>0$ is chosen sufficiently small. Inequality (\ref{supAL2}) then follows by applying
the dominated convergence theorem in (\ref{596}).

The proof of (\ref{supAL}) for any $\rho>0$ and the uniform convergence in (\ref{supAL}) for $\rho \in [\rho_0,1]$ for any
$\rho_0\in (0,1]$ follows by similar arguments. It suffices to observe that $\bar\beta' \leq 1-C\sqrt{\bar\beta -1}$ for
some $C>0$ uniformly in $\rho\in [\rho_0, 1]$ and $\bar\beta>1$ sufficiently small. This concludes the proof of Theorem \ref{T:cptst}.
\qed
\medskip

\noindent
{\bf Remark.} Note that the change of measure argument here applies equally well to any random walks $X$ and $Y$
with an identical symmetric transition kernel.

\subsection{Proof of Theorem \ref{T:cptst} in continuous time: $d=4$}
As in $d\geq 5$, proving Theorem \ref{T:cptst} reduces to proving $\varrho<1$ (see (\ref{rhocts})) for
appropriate choices of $\bar\beta$ and $L$ depending on the diffusion constant $\rho$. Since
$\E^Y_0[\bar Z^{\bar\beta}_{t,Y}]$ is the partition function of a homogeneous pinning model with parameter
$\bar\beta\geq 1$ and return time distribution $K_{1+\rho}(t)\sim C t^{-2}$, by comparing $K_{1+\rho}$ with a return
time distribution $K'$ which is stochastically smaller than $K_{1+\rho}$ and has finite first moment, we see that
(\ref{517}) also holds in $d=4$. Therefore setting $L=(\bar\beta-1)^{-1}$ as in $d\geq 5$, we have
$\sup_{0\leq t\leq L} \bar A^{\bar\beta, \rm pin2}_t \leq C<\infty$, and analogous to (\ref{rhodec4}), we have
\be
\varrho \leq \sum_{i=0}^{L^{1-\eps}} \frac{C}{(L-i)^{2\gamma-1}}
+ \sum_{i=L^{1-\eps}}^{L-1} \frac{C \bar A^{\bar\beta, \rm pin2}_i}{(L-i)^{2\gamma-1}}
\leq \frac{CL^{1-\eps}}{(L-L^{1-\eps})^{2\gamma-1}} + C L^{2-2\gamma} \!\!\!\!\!\!
\sup_{L^{1-\eps}\leq t\leq L}\bar A^{\bar\beta, \rm pin2}_t,
\ee
where $\eps>0$, $\gamma \in (0,1)$ is chosen so that $2\gamma -1>1-\eps$, and $C\in (0,\infty)$ is independent
of $\bar\beta \in [1,2]$ and is furthermore uniformly bounded for $\rho \in [0,1]$. Therefore, to show
$\beta^*_{\rm c} >\beta^{\rm ann}_{\rm c}$ for any $\rho>0$, it suffices to show
\be\label{529}
\lim_{\bar\beta \downarrow 1} L^{2-2\gamma}\!\!\!\! \sup_{L^{1-\eps\leq t\leq L}} \bar A^{\bar\beta, \rm pin2}_t = 0.
\ee
On the other hand, to show that for any $\delta>0$, there exists $a_\delta>0$ such that
\be\label{530}
\beta^*_{\rm c} - \beta^{\rm ann}_{\rm c} \geq a_\delta \rho^{1+\delta} \qquad \forall\ \rho \in [0,1],
\ee
it suffices to show that: (1) the convergence in (\ref{529}) is uniform for $\rho \in [\rho_0,1]$ for any
$0<\rho_0\leq 1$, which implies that $\inf_{\rho \in [\rho_0,1]}(\beta^*_{\rm c}-\beta^{\rm ann}_{\rm c})>0$;
(2) for $\bar\beta = 1+ \rho^{1+\delta}$ and $L=(\bar\beta-1)^{-1} = \rho^{-1-\delta}$,
\be\label{531}
\lim_{\rho\downarrow 0} L^{2-2\gamma}\!\!\!\! \sup_{L^{1-\eps}\leq t\leq L} \bar A^{\bar\beta, \rm pin2}_t = 0,
\ee
which implies that for some $\rho_0 \in (0,1]$, $\bar\beta^*_{\rm c} -1= G_{1+\rho}(\beta^*_{\rm c}-\beta^{\rm ann}_{\rm
  c})\geq \rho^{1+\delta}$ for all $\rho \in [0, \rho_0]$.

Proceeding exactly as in the $d\geq 5$ case, we note that (\ref{168}) still holds in $d= 4$. By the choice
$h= \frac{\sqrt \rho}{\sqrt L} = \rho^{1+\delta/2}$, there exists $\rho_1 \in (0,1)$ such that
\be
\bar\beta' = \frac{(1+\rho)\bar\beta}{1+\rho+h} = \frac{(1+\rho)(1+\rho^{1+\delta})}{1+\rho+\rho^{1+\delta/2}}\leq
1-\rho^{1+\delta/2}/2 \leq e^{-\rho^{1+\delta/2}/2}\qquad \forall \ \ \rho \in [0, \rho_1].
\ee
If we choose $K_*$ to be a return time distribution with $\int_0^\infty K_*(s)ds=1$
and $K_*(s) \sim C s^{-2}$ such that $K_*$ stochastically dominates $K_{1+\rho+h}$ for all $\rho, h\in [0,1]$, and
let $\iota^*$ be a renewal process on $[0,\infty)$ with return time distribution $K_*$, then
\begin{eqnarray}
0\leq \lim_{\rho\downarrow 0} L^{2-2\gamma}\!\!\!\! \sup_{L^{1-\eps}\leq t\leq L} \bar A^{\bar\beta, \rm pin2}_t
&\leq& C \lim_{\rho\downarrow 0} \rho^{-(1+\delta)(2-2\gamma)}
\E^{K_*}\left[\exp\left\{-\frac{1}{2}\rho^{1+\delta/2}\, \big|\iota^* \cap
    [0,\rho^{-(1+\delta)(1-\eps)}]\big|\right\}\right]^\gamma \nonumber \\
&=& C \lim_{M\to\infty} M^{\frac{2-2\gamma}{1-\eps}}
\E^{K_*}\left[\exp\left\{-\frac{1}{2}M^{-\frac{1+\delta/2}{(1+\delta)(1-\eps)}} \big|\iota^* \cap [0,M]|
  \right\}\right]^\gamma =0 \nonumber
\end{eqnarray}
where we applied Proposition \ref{P:rewest} with $\delta_1 = \frac{1+\delta/2}{(1+\delta)(1-\eps)}$ and $1-\delta_2 =
\frac{2-2\gamma}{\gamma(1-\eps)}$, which satisfy the condition $0<\delta_1<\delta_2<1$ if $\eps>0$ is small and
$\gamma$ is then chosen sufficiently close to 1. This proves (\ref{531}).

The proof of (\ref{529}) for any $\rho>0$ and the uniform convergence therein for $\rho \in [\rho_0,1]$ for any
$\rho_0\in (0,1]$ follows by similar arguments. It suffices to note that for each $\rho>0$, there exists $C>0$
and $\bar\beta_0 >1$ such that $\bar\beta' \leq 1-C\sqrt{\bar\beta -1}$ for all $\bar\beta \in [1, \bar\beta_0]$.
Furthermore, $C$ and $\bar\beta_0$ can be chosen uniformly for $\rho\in [\rho_0, 1]$ for any $\rho_0>0$. The
rest of the proof proceeds exactly as for $d=4$ in the discrete time case.
\qed

\appendix
\section{A renewal process estimate}
The following proposition complements Proposition A.2 in \cite{DGLT07}  for the case $\alpha=1$.

\bp\label{P:rewest}
Let $\iota^*=\{\iota_0=0, \iota_1, \cdots\}$ be a renewal process on $\N_0$ with inter-arrival probability distribution
$K_*$ satisfying $\sum_{n\in\N} K_*(n)=1$ and $K_*(n) \sim C n^{-2}$ as $n\to\infty$. Then for any $c>0$ and
$0<\delta_1 <\delta_2<1$, we have
\be\label{A1}
\lim_{N\to\infty} N^{1-\delta_2} \E^{K_*}\left[\exp\left\{-c N^{-\delta_1} \big|\iota^*\cap [0,N]\big| \right\}\right] = 0.
\ee
The same result holds if $\iota^*$ is a renewal process on $[0,\infty)$ with inter-arrival distribution $K_*$ satisfying
$\int_0^\infty K_*(s)ds=1$ and $K_*(s)\sim C s^{-2}$ as $s\to\infty$.
\ep
{\bf Proof.} Let $\delta_3 \in (\delta_1, \delta_2)$. Note that
\be\label{A2}
\E^{K_*}\left[\exp\left\{-c N^{-\delta_1} \big|\iota^*\cap [0,N]\big| \right\}\right]
\leq \P\left(0\leq |\iota^*\cap [0,N]| < N^{\delta_3}\right) + e^{-c N^{\delta_3-\delta_1}}.
\ee
Let $(U_i)_{i\in\N}$ be i.i.d.\ random variables with distribution $K_*$. By our assumption on $K_*$, for each $\alpha \in (0,1)$,
we can find a constant $C_\alpha>0$ and i.i.d.\ stable subordinators $(V_i)_{i\in\N}$ with exponent $\alpha$, i.e., $\P(V_1>0)=1$ and
$V_1 \stackrel{\rm law}{=} \sum_{i=1}^n V_i / n^{1/\alpha}$, such that $\P(U_1>s)\leq \P(V_1+C_\alpha>s)$ for all $s>0$. Therefore,
for $\alpha \in (\delta_3, 1)$,

\begin{eqnarray}
&& \P\left(0\leq |\iota^*\cap [0,N]| < N^{\delta_3}\right) = \P\left(\sum_{n=1}^{N^{\delta_3}} U_n > N\right) \nonumber \\
&\leq& \P\left(\sum_{n=1}^{N^{\delta_3}} (V_n+C_\alpha) >N \right) = \P\left(\sum_{n=1}^{N^{\delta_3}} V_n > N-C_\alpha N^{\delta_3}\right)
=\P\left(V_1 > N^{1-\delta_3/\alpha} - C_\alpha N^{\delta_3(1-1/\alpha)} \right) \nonumber \\
&\leq& CN^{\delta_3-\alpha}, \label{A3}
\end{eqnarray}
where we used the fact that $\P(V_1>x) \sim C x^{-\alpha}$ as $x\to\infty$. It is easy to see that (\ref{A1}) follows from
(\ref{A2}) and (\ref{A3}) if we choose $\alpha\in (0,1)$ such that $1-\delta_2+\delta_3-\alpha<0$. The case when $\iota^*$ is
a renewal process on $[0,\infty)$ can be treated identically.
\qed
\bigskip

\noindent
{\bf Acknowledgment} R.S.\ thanks J\"urgen G\"artner and Thomas Mountford for helpful discussions on the parabolic
Anderson model with moving catalysts. We thank F.L.Toninelli for interesting discussions, and the referees for helpful
comments and suggestions. R.S.\ was supported by a postdoctoral position in the DFG Forschergruppe 718 {\it Analysis and
Stochastics in Complex Physical Systems}.

\end{document}